\documentclass[5p]{elsarticle}

\usepackage[intlimits]{amsmath}
\usepackage{amsthm,amssymb,amsfonts}
\usepackage[english]{babel}
\usepackage{xcolor}

\usepackage{graphics} 
\usepackage{epsfig} 
\theoremstyle{plain}
\newtheorem{assumption}{Assumption}

\newtheorem{lemma}{Lemma}
\newtheorem{theorem}{Theorem}
\newtheorem{proposition}{Proposition}
\theoremstyle{remark}
\newtheorem{remark}{Remark}
\biboptions{sort&compress}
\DeclareMathOperator*{\esssup}{ess\,sup}

\journal{Automatica}

\begin{document}
	
\begin{frontmatter}
		
\title{Observer-Based Output-Feedback Backstepping Stabilization of Continua of Hyperbolic 
PDEs and Application to Large-Scale $n+m$ Coupled Hyperbolic PDEs\tnoteref{t1}} 
\tnotetext[t1]{Funded by the European Union (ERC, C-NORA, 101088147). Views and 
opinions expressed are however those of the authors only and do not necessarily reflect those of 
the European Union or the European Research Council Executive Agency. Neither the European 
Union nor the granting authority can be held responsible for them.}
		
\author[1]{Jukka-Pekka Humaloja}
\author[1]{Nikolaos Bekiaris-Liberis}
		
\affiliation[1]{organization={Department of Electrical and Computer Engineering, Technical 
University of Crete},
	addressline={University Campus, Akrotiri}, 
	city={Chania},
	postcode={73100}, 
	country={Greece}}
		
\begin{abstract}
We develop a non-collocated, observer-based output-feedback law for a class of continua of 
linear hyperbolic PDE systems, which are viewed as the continuum version of $n+m$, general 
heterodirectional hyperbolic systems as $n\to\infty$. The design relies on the introduction of a 
novel, continuum PDE backstepping transformation, which enables the construction of a 
Lyapunov functional for the estimation error system. Stability under the observer-based 
output-feedback law is established by using the Lyapunov functional construction for the 
estimation error system and proving well-posedness of the complete closed-loop system, which 
allows utilization of the separation principle. 

Motivated by the fact that the continuum-based 
designs may provide computationally tractable control laws for large-scale, $n+m$ systems, we 
then utilize the control/observer kernels and the observer constructed for the continuum system 
to introduce an output-feedback control design for the original $n+m$ system. We establish 
exponential stability of the resulting closed-loop system, which consists of a mixed 
$n+m$-continuum PDE system (comprising the plant-observer dynamics), introducing a virtual 
continuum system with resets, which 
enables utilization of the continuum approximation property of the solutions of the $n+m$ system 
by its continuum counterpart (for large $n$). We illustrate the potential computational 
complexity/flexibility benefits of our approach via a numerical example of stabilization of a 
large-scale $n+m$ system, for which we employ the continuum observer-based controller, while 
the continuum-based stabilizing control/observer kernels can be computed in closed form.
\end{abstract}
		
		
\end{frontmatter}

\section{Introduction}

\subsection{Motivation}

Continua of hyperbolic PDE systems can be viewed as continuum versions of certain, large-scale 
hyperbolic systems, featuring a large number of PDE state components  \cite{HumBek24arxivc, 
	HumBek25b}. A specific, theoretically and practically significant case of the latter, is the 
class  of large-scale, $n+m$, heterodirectional, linear hyperbolic systems, which may be utilized 
to describe, for example, the dynamics of blood \cite{BikPhd, ReyMer09}, traffic 
\cite{ZhaLua22, HerKla03, YuHKrs21}, and water \cite{BasCorBook,DiaDia17} flow 
networks,  as well as the dynamics of epidemics transport \cite{GuaPri20, KitBes22}. For example, 
PDE-based traffic flow models for multi-lane traffic give rise to $n+m$ systems with $(n-m)/2$ 
lanes in free-flow and $m$ lanes in congested conditions, respectively \cite{HerKla03}. Whereas 
PDE models of arterial networks describing blood flow from the heart all the way through to a 
brachial artery (where a measurement can be obtained in a non-invasive manner) consist of 
interconnected, $n+m$ hyperbolic systems \cite{ReyMer09}.

In fact, certain control  designs developed for stabilization of continua of hyperbolic PDE 
systems can be utilized for  stabilization of the corresponding large-scale system when the 
number of state components is  sufficiently large \cite{HumBek24arxivc, HumBek25b}. This is 
an important feature as it  allows construction of stabilizing control kernels whose computational 
complexity does not grow  with the number of state components \cite{HumBek25, 
HumBek25b}. A natural next step  is to introduce a dual approach in which one constructs 
observers and observer kernels for continua of hyperbolic 
systems, which could, in principle, be utilized as (approximate) observers and observer kernels for 
the large-scale system counterpart. This key idea provides to the designer the flexibility to 
construct and compute both observer kernels and observer dynamics, essentially, 
independently of $n$, which has the potential to achieve design of computationally tractable 
observer-based control laws for large-scale PDE systems. Motivated by this and the fact that 
neither an observer-based output-feedback design for continua systems of hyperbolic PDEs is 
available nor it has been applied to control of large-scale $n+m$ systems, in the present paper we 
address the problems of observer and output-feedback designs for such a class of continua of 
hyperbolic systems, as well as their application to large-scale $n+m$ systems.

\subsection{Literature}

Full-state feedback laws for a class of continua of linear hyperbolic PDEs have been recently 
developed in \cite{AllKrs25} and  \cite{HumBek24arxivc}. In particular, \cite{AllKrs25} first 
addressed the problem of stabilization of a continuum version of the $n+1$ systems considered in 
\cite{DiMVaz13} as $n\to\infty$, whereas in  \cite{HumBek24arxivc} we developed a feedback 
control design approach for stabilization of the continuum counterpart (as $n\to\infty$) of the 
$n+m$ hyperbolic systems considered in \cite{HuLDiM16}. The fact that the control design 
procedure developed for the continuum system in \cite{AllKrs25} can be employed for 
stabilization of the large-scale $n+1$ (for finite $n$) system, as it may provide stabilizing control 
kernels that can approximate to arbitrary accuracy the exact (constructed directly for the 
large-scale $n+1$ system) backstepping kernels (for sufficiently large $n$), was established in 
\cite{HumBek25b, HumBek25b}. There exists no result heretofore addressing the problems 
of observer and output-feedback designs for such classes of continua of hyperbolic PDE systems. 

The backstepping-based output-feedback stabilization problem of $n+m$ hyperbolic systems 
has been solved  in
\cite{HuLDiM16}, where the proposed control law involves solving the $n+m$ control and 
observer kernel equations, and constructing a Luenberger-type observer for the $n+m$ system. 
More recently,  control design methods have been developed for other types of $n+m$ systems 
in \cite{AurBri24, Aur24, EndGab24}. Moreover, state-feedback stabilization of various types of 
$n+m$ systems has been considered, e.g., in \cite{AurBre22, AurDiM16, CorHuL17, DiMArg18, 
HuLVaz19}. Such approaches may result in high derivation complexity and computational burden 
of the respective control laws, as these increase with $n$ and $m$ (see, e.g., \cite{Aur24}), due 
to the requirements of obtaining the solutions to the control/observer kernels, as well as of 
implementing the observer dynamics. To overcome these potential computational obstacles, 
methods based 
on, e.g., neural networks \cite{BhaShi24, WanDia24arxiv, QiJZha24} and late-lumping 
implementations \cite{AurMor19} have been considered 
to efficiently approximate a given control design. Our approach can be viewed as an 
alternative to the above, in that the main goal is to avoid the increase in computational complexity 
with respect to $n$, as well as to provide the user implementation flexibility, via continuum 
approximation of the 
backstepping kernels and the observer dynamics, which is a problem heretofore unattempted.

\subsection{Contributions}

In the present paper we develop a backstepping-based observer design methodology for a class 
of continua systems, which are viewed as the continuum version of $n+m$, linear hyperbolic 
systems as $n\to\infty$. Specifically, we address the dual to the control design problem from  
\cite{HumBek24arxivc} in which we consider availability of $m$ measurements, anti-collocated to 
the boundary where control is applied. Introducing a suitable target system we derive the 
continuum kernel equations, which are shown to be well-posed by recasting them in the form of 
the control kernel equations from \cite{HumBek24arxivc}. Our choice of the target system 
enables construction of a Lyapunov functional for the estimation error system. We then introduce 
the respective observer-based output-feedback design combining the observer design developed 
here with the control design developed in  \cite{HumBek24arxivc}. The key in establishing stability 
of the complete closed-loop system is to show its well-posedness, which in turn allows us to 
employ the separation principle. The well-posedness proof relies on recasting the complete 
closed-loop system as an abstract system in output-feedback form and deriving the respective 
transfer function matrix. 

We then utilize as basis the observer and control designs developed for the continuum to 
introduce an observer-based control design methodology for the large-scale, $n+m$ system 
counterpart. The methodology consists of two main ingredients, namely, the construction of a 
continuum-based observer (that is employed for estimation of the $n+m$ system's state) and the 
construction of continuum-based observer kernels. The key element in the methodology lies in 
the execution of the first step that constitutes the method introduced here 
different from the respective method in \cite{HumBek24arxivc}, which considers the full-state 
feedback control design case. In particular, we establish that an observer constructed on the 
basis of the continuum system can provide accurate estimates of the state of the $n+m$ system 
provided that $n$ is large and that a suitable sampling is applied. Moreover, we study the stability 
properties of the resulting, unexplored closed-loop system, consisting of both an $n+m$ and a 
continuum PDE system, introducing a virtual continuum system with resets at properly chosen 
time instants. This 
enables utilization of the continuum approximation property of the solutions of the $n+m$ system 
by its continuum counterpart. Specifically, we show that the (augmented by the virtual continuum) 
closed-loop system can be expressed as consisting of an exponentially stable nominal part that is 
affected by an additive (state-dependent) perturbation, which does not destroy exponential 
stability when $n$ is sufficiently large, as this allows, at each resetting time instant, to reset the 
state of the augmented system to a state of smaller magnitude. We also present a numerical 
example with consistent simulation results, in which the continuum-based control/observer 
kernels can be computed in closed form and the continuum observer is employed, illustrating the 
potential benefits of our approach in computational complexity/flexibility.

\subsection{Organization}

In Section~\ref{sec:infm}, we present the class of continua of hyperbolic PDEs considered and 
the proposed observer-based output-feedback law. In Section~\ref{sec:obskappr}, we employ 
continuum-based, backstepping control/observer kernels in designing an observer-based control 
law for $n+m$ 
hyperbolic systems. In Section~\ref{sec:obsappr}, we propose a continuum observer-based 
control law for $n+m$ hyperbolic PDEs, where both the observer design and backstepping 
kernels are based 
on the continuum approximation. In Section~\ref{sec:ex}, we present a numerical example and 
consistent simulation results. In Section~\ref{sec:conc}, we provide concluding remarks and 
discuss potential topics of future research.

\subsection{Notation}

We use the standard notation $L^2(\Omega; \mathbb{R})$ for real-valued
Lebesgue integrable functions on any domain $\Omega \subset \mathbb{R}^d$ for some 
$d \geq 1$. The notations $L^\infty(\Omega;\mathbb{R}),C(\Omega; \mathbb{R})$, and 
$C^1(\Omega; \mathbb{R})$ denote essentially bounded, continuous, and
continuously differentiable functions, respectively, on $\Omega$. We 
denote vectors and matrices by bold symbols, and any $n,m \in 
\mathbb{N}$, we denote by $E$ the 
space $L^2([0,1]; \mathbb{R}^{n+m})$ equipped with the inner product 
\begin{align}
\label{eq:eip}
\left\langle \left( 
\begin{smallmatrix}
\mathbf{u}_1 \\ \mathbf{v}_1
\end{smallmatrix}
\right), \left( 
\begin{smallmatrix}
\mathbf{u}_2\\ \mathbf{v}_2
\end{smallmatrix}
\right) \right\rangle_E
& = \nonumber \\ \int\limits_0^1 \left(\frac{1}{n}
\sum_{i=1}^nu_1^i(x)u_2^i(x) +
  \sum_{j=1}^mv_1^j(x)v_2^j(x)\right)dx,
&
\end{align}
which induces the norm $\left\|
\cdot
\right\|_E = \sqrt{\left<\cdot,\cdot\right>_E}$. We also define the continuum version of $E$ as 
$n\to\infty$ by
$E_c = L^2([0,1];L^2([0,1];\mathbb{R}))\times L^2([0,1];\mathbb{R}^m)$, (i.e., $\mathbb{R}^n$ 
becomes $L^2([0,1];\mathbb{R})$ as $n\to\infty$) equipped with the inner 
product
\begin{align}
\label{eq:ecip}
\left\langle \left( 
\begin{smallmatrix}
u_1 \\ \mathbf{v}_1
\end{smallmatrix}
\right), \left( 
\begin{smallmatrix}
u_2\\ \mathbf{v}_2
\end{smallmatrix}
\right) \right\rangle_{E_c}
& = \nonumber \\
 \int\limits_0^1\left(\int\limits_0^1 u_1(x,y)u_2(x,y)dy + \sum_{j=1}^mv_1^j(x)v_2^j(x)\right)dx, &
\end{align}
which coincides with $L^2([0,1]^2;\mathbb{R})\times L^2([0,1];\mathbb{R}^m)$.
Moreover, the transform $\mathcal{F} =
\operatorname{diag}(\mathcal{F}_{n}, I_m)$ is an isometry from $E$ to $E_c$, where 
$\mathcal{F}_n$ maps any $\mathbf{b} \in \mathrm{R}^n$ to $L^2([0,1]; \mathbb{R})$ as
\begin{equation}
	\label{eq:Fn}
	\mathcal{F}_n\mathbf{b} =\sum_{i=1}^n b_i \chi_{((i-1)/n,i/n]},
\end{equation}
where $\chi_{((i-1)/n,i/n]}$ denotes the indicator function of the interval $((i-1)/n,i/n]$. The adjoint 
of $\mathcal{F}$ is of the form $\mathcal{F}^* = \operatorname{diag}(\mathcal{F}^*_n,I_m)$ with
\begin{equation}
	\label{eq:Fns}
	\mathcal{F}_n^{*}g = \left( n\int\limits_{(i-1)/n}^{i/n}g(y)dy \right)_{i=1}^n,
\end{equation}
where the $i$-th component is the mean value of the function $g \in L^2([0,1]; \mathbb{R})$ over 
the interval $[(i-1)/n,i/n]$. Finally, we denote by $\mathcal{L}(X,Z)$ the space of bounded linear 
operators from any normed space $X$ to any normed space $Z$, and 
$\|\cdot\|_{\mathcal{L}(X,Z)}$ denotes the corresponding operator norm.

We say that a system is exponentially stable on a normed space $Z$ if, for any initial condition 
$z_0\in Z$, the (weak) solution 
$z \in C([0,\infty); Z)$ of the system satisfies $\|z(t)\|_Z \leq 
Me^{-ct}\|z_0\|_E$ for some constants $M,c > 0$ that are independent of $z_0$. Finally, we 
denote by $\mathcal{T}$ the triangular set
\begin{align}
\mathcal{T} & = \left\{ (x,\xi) \in [0,1]^2: 0 \leq \xi \leq x \leq 1 \right\}.
\end{align}

\section{Observer-Based Output-Feedback Stabilization of Continua Systems of Hyperbolic 
PDEs} \label{sec:infm}

\subsection{Continua Systems of Hyperbolic PDEs}

The considered class of continuum systems is of the form
\begin{subequations}
	\label{eq:infm}%
	\begin{align}
		u_t(t,x,y) + \lambda(x,y)u_x(t,x,y)   & = \nonumber \\
		\int\limits_0^1\sigma(x,y,\eta)u(t,x,\eta)d\eta +
		\mathbf{W}(x,y)\mathbf{v}(t,x), & \label{eq:infm1} \\
		\mathbf{v}_t(t,x) - \pmb{\Lambda}_-(x)\mathbf{v}_x(t,x)
		& =  \nonumber \\
		\int\limits_0^1\pmb{\Theta}(x,y)u(t,x,y)dy + \pmb{\Psi}(x)\mathbf{v}(t,x),
		& \label{eq:infm2}
	\end{align}
\end{subequations}
with boundary conditions 
\begin{subequations}
	\label{eq:infmbc}%
	\begin{align}
		u(t,0,y) & = \mathbf{Q}(y)\mathbf{v}(t,0), \\
		\mathbf{v}(t,1) & = \int\limits_0^1 \mathbf{R}(y)u(t,1,y)dy + \mathbf{U}(t),
	\end{align}
\end{subequations}
and output $\mathbf{Y}(t) = \mathbf{v}(t,0)$, for almost every $y \in [0,1]$. Here we employ the 
matrix notation for $\mathbf{v}, \mathbf{U}, \mathbf{Y}, \pmb{\Lambda}_-, \pmb{\Theta},
\pmb{\Psi}, \mathbf{W},\mathbf{Q}$, and $\mathbf{R}$ for the sake of 
conciseness, that is, $\mathbf{v} = \left( v^j \right)_{j=1}^m$,
$\mathbf{U} = \left( U^j \right)_{j=1}^m$, $\mathbf{Y} = \left( Y^j \right)_{j=1}^m$, and the 
parameters are as follows.
\begin{assumption}
	\label{ass:infm}
	The parameters of \eqref{eq:infm}, \eqref{eq:infmbc} are such that
	\begin{subequations}
		\label{eq:infmparam}
		\begin{align}
			\pmb{\Lambda}_- & = \operatorname{diag}(\mu_j)_{j=1}^m
			\in C^1([0,1]; \mathbb{R}^{m\times m}), \\
			\pmb{\Theta} & = \left( \Theta_j \right)_{j=1}^m \in
			C([0,1];L^2([0,1];\mathbb{R}^m)), \\
			\pmb{\Psi} & = \left( \Psi_{i,j} \right)_{i,j=1}^m \in C([0,1];
			\mathbb{R}^{m\times m}), \\
			\mathbf{W} & = 
			\begin{bmatrix}
				W_1 & \cdots  & W_m
			\end{bmatrix} \in C([0,1];L^2([0,1];\mathbb{R}^{1\times m})), \\
			\mathbf{Q} & = 
			\begin{bmatrix}
				Q_1 & \cdots & Q_{m}
			\end{bmatrix} \in L^2([0,1]; \mathbb{R}^{1\times m}), \\
			\mathbf{R} & =  \left( R_j \right)_{j=1}^m \in L^2([0,1];\mathbb{R}^m),
		\end{align}
	\end{subequations}
	with $\lambda  \in C^1([0,1]^2;\mathbb{R})$ and  $\sigma \in
	C([0, 1]; L^2([0,1]^2;\mathbb{R}))$. Moreover, $\mu_m(x) > 0$ and $\lambda(x,y) > 0$ 
	for all $x,y \in [0,1]$, and additionally
	\begin{equation}
		\label{eq:muass}
		\min_{x\in [0,1]}\mu_j(x) > \max_{x\in[0,1]}\mu_{j+1}(x),
	\end{equation}
	for all $j = 1,\ldots, m-1$. Finally, $\psi_{j,j} = 0$ for all $j = 
	1,\ldots,m$.\footnote{This comes without loss of generality, as such terms
		can be removed using a change of variables (see, e.g., \cite{HuLDiM16, HuLVaz19}).}
\end{assumption}

\subsection{Control Law and Observer Design}

The control law to stabilize \eqref{eq:infm}, \eqref{eq:infmbc} is of the form 
\begin{align}
	\label{eq:infmUobs}
	\mathbf{U}(t) &  =  \int\limits_0^1\int\limits_0^1
	\mathbf{K}(1,\xi,y)\hat{u}(t,\xi,y)dyd\xi \nonumber \\
	& \qquad + \int\limits_0^1 \mathbf{L}(1,\xi)\hat{\mathbf{v}}(t,\xi) d\xi 
	-\int\limits_0^1\mathbf{R}(y)\hat{u}(t,1,y)dy,
\end{align}
where $\mathbf{K}, \mathbf{L}$ (satisfying \eqref{eq:infmk}--\eqref{eq:infmkbca1}) are the 
backstepping control kernels recalled in \ref{app:sf} and 
$\hat{u}, \hat{\mathbf{v}}$ is the observer state satisfying the 
dynamics
\begin{subequations}
\label{eq:infmobs}%
\begin{align}
	\resizebox{!}{.028\columnwidth}{$	
	\hat{u}_t(t,x,y) + \lambda(x,y)\hat{u}_x(t,x,y) + \mathbf{P}_+(x,y)(\hat{\mathbf{v}}(t,0) - 
	\mathbf{v}(t,0))$} & = 
	\nonumber \\
	\int\limits_0^1\sigma(x,y,\eta)\hat{u}(t,x,\eta)d\eta +
	\mathbf{W}(x,y)\hat{\mathbf{v}}(t,x), & \\
	\hat{\mathbf{v}}_t(t,x) - \pmb{\Lambda}_-(x)\hat{\mathbf{v}}_x(t,x) + 
	\mathbf{P}_-(x)(\hat{\mathbf{v}}(t,0) - 
	\mathbf{v}(t,0))
	& =  \nonumber \\
	\int\limits_0^1\pmb{\Theta}(x,y)\hat{u}(t,x,y)dy + \pmb{\Psi}(x)\hat{\mathbf{v}}(t,x),
	&
\end{align}
\end{subequations}
with boundary conditions 
\begin{subequations}
	\label{eq:infmobsbc}%
	\begin{align}
		\hat{u}(t,0,y) & = \mathbf{Q}(y)\mathbf{v}(t,0), \\
		\hat{\mathbf{v}}(t,1) & = \int\limits_0^1 \mathbf{R}(y)\hat{u}(t,1,y)dy + \mathbf{U}(t),
	\end{align}
\end{subequations}
where $\mathbf{P}_+, 
\mathbf{P}_-$ are given by
\begin{subequations}
	\label{eq:infmP}
	\begin{align}
		\mathbf{P}_+(x,y) & = \mathbf{M}(x,0,y)\pmb{\Lambda}_-(0), \\
		\mathbf{P}_-(x) & = \mathbf{N}(x,0)\pmb{\Lambda}_-(0),
	\end{align}
\end{subequations}
where $\mathbf{M} \in L^\infty(\mathcal{T};L^2([0,1];\mathbb{R}^{1\times m})), \mathbf{N}  \in 
L^\infty(\mathcal{T};\mathbb{R}^{m\times m}))$, is the solution to the observer kernel equations
\begin{subequations}
	\label{eq:infmobsk}
	\begin{align}
		\lambda(x,y)\mathbf{M}_x(x,\xi,y) - \mathbf{M}_\xi(x,\xi,y)\pmb{\Lambda}_-(\xi) -  
		\mathbf{M}(x,\xi,y)\pmb{\Lambda}_-'(\xi) & = \nonumber \\	
		\int\limits_0^1 \sigma(\xi,y,\eta)\mathbf{M}(x,\xi,\eta)d\eta + 
		\mathbf{W}(\xi,y)\mathbf{N}(x,\xi), \\
		\pmb{\Lambda}_-(x)\mathbf{N}_x(x,\xi) + \mathbf{N}_\xi(x,\xi)\pmb{\Lambda}_-(\xi) + 
		\mathbf{N}(x,\xi)\pmb{\Lambda}_-'(\xi) & = \nonumber \\
		\int\limits_0^1 \pmb{\Theta}(\xi,y)\mathbf{M}(x,\xi,y)dy + \pmb{\Psi}(\xi)·\mathbf{N}(x,\xi),
	\end{align}
\end{subequations}
with boundary conditions
\begin{subequations}
	\label{eq:infmobskbc}
	\begin{align}
		\mathbf{W}(x,y) & = \mathbf{M}(x,x,y)\pmb{\Lambda}_-(x) + \lambda(x,y)\mathbf{M}(x,x,y),
		\label{eq:infmobskbc1} \\
		\pmb{\Psi}(x) & = \mathbf{N}(x,x)\pmb{\Lambda}_-(x) - \pmb{\Lambda}_-(x)\mathbf{N}(x,x),
		\label{eq:infmobskbc2} \\
		N_{i,j}(1,\xi) & =  \int\limits_0^1 R_i(y)M_j(1,\xi,y)dy, \quad  \forall i\geq j, 
		\label{eq:infmobskbc3} \\
		N_{i,j}(x,0) & = n_{i,j}(x), \quad \forall i < j, \label{eq:infmobskbc4}
	\end{align}
\end{subequations}
where $n_{i,j}(x)$ are arbitrary due to \eqref{eq:infmobskbc4} being an artificial boundary 
condition to guarantee well-posedness of the observer kernel equations; similarly to the $n+m$ 
case in \cite{HuLDiM16}.  We choose $n_{i,j}(x)$ 
such that 
\begin{equation}
	\label{eq:nij}
	n_{i,j}(0) = \frac{\Psi_{i,j}(0)}{\mu_j(0) - \mu_i(0)},
\end{equation}
in order to make the artificial boundary condition compatible with 	\eqref{eq:infmobskbc2} at 
$(0,0)$. Note that the boundary conditions for $\mathbf{N}(1,1)$ are, in general, overdetermined 
due to \eqref{eq:infmobskbc2} and \eqref{eq:infmobskbc3},  \eqref{eq:infmobskbc1}, which stems 
potential discontinuities in the $\mathbf{N}$ kernels.

\subsection{Well-Posedness of Observer Kernel Equations}

\begin{theorem}
	\label{thm:obskwp}
	Under Assumption~\ref{ass:infm}, equations \eqref{eq:infmobsk}--\eqref{eq:nij} have a 
	well-posed solution $\mathbf{M} \in
	L^\infty(\mathcal{T};L^2([0,1];\mathbb{R}^{1\times m}))$, $\mathbf{N}  \in 
	L^\infty(\mathcal{T};\mathbb{R}^{m\times m}))$. Moreover, the solution is piecewise 
	continuous 
	in $(x,\xi) \in \mathcal{T}$, where the set of discontinuities is of measure zero.
	\begin{proof}
		We first transform the kernel equations \eqref{eq:infmobsk}--\eqref{eq:nij} into an 
		analogous form with the control kernel equations, which have been shown to be well-posed 
		in 
		\cite[Thm 2]{HumBek24arxivc}. This transformation is achieved
		by introducing alternative variables 
		\begin{subequations}
			\label{eq:MNalt}
			\begin{align}
				\bar{\mathbf{M}}(\chi,\zeta,y) & = \mathbf{M}(1-\zeta,1-\chi,y) = \mathbf{M}(x,\xi,y), \\
				\bar{\mathbf{N}}(\chi,\zeta) & = \mathbf{N}(1-\zeta,1-\chi) = \mathbf{N}(x,\xi),
			\end{align}
		\end{subequations}
		and analogously introducing the alternative parameters $\bar{\mu}, \bar{\lambda}, 
		\bar{\sigma},\bar{\Theta}$, 
		$\bar{W}, \bar{\Psi}, \bar{n}_{i,j}$ based on the coordinate transform $(x,\xi) \to (1-\zeta, 
		1-\chi)$, 
		so that $(\chi,\zeta)$ is the mirror image of $(x,\xi)$ with respect to the line $x+\xi = 1$, and 
		hence, $(x,\xi) \in \mathcal{T}$ is mirrored into $(\chi,\zeta) \in \mathcal{T}$. Inserting these 
		new 
		coordinates into the observer kernel equations \eqref{eq:infmobsk}, \eqref{eq:infmobskbc}  
		and 
		noting $\partial_x = -\partial_\zeta, \partial_\xi = -\partial_\chi$, we get the componentwise 
		equations, for $i,j=1,\ldots,m$,
		\begin{subequations}
			\label{eq:infmobskalt}
			\begin{align}
				\resizebox{.93\columnwidth}{!}{$ \displaystyle 
					\bar{\mu}_j(\chi)\partial_\chi \bar{M}_j(\chi,\zeta,y) - 
					\bar{\lambda}(\zeta,y)\partial_\zeta 
					\bar{M}_j(\chi,\zeta,y) + 
					\bar{M}_j(\chi,\zeta,y)\bar{\mu}_j'(\chi)$} & = \nonumber \\
				\int\limits_0^1 \bar{\sigma}(\chi,y,\eta)\bar{M}_j(\chi,\zeta,\eta)d\eta + 
				\sum_{\ell=1}^{m}\bar{W}_\ell(\chi,y)\bar{N}_{\ell,j}(\chi,\zeta), \\
				\resizebox{.93\columnwidth}{!}{$ \displaystyle 
					\bar{\mu}_j(\chi)\partial_\chi \bar{N}_{i,j}(\chi,\zeta) + \bar{\mu}_i(\zeta)\partial_\zeta 
					\bar{N}_{i,j}(\chi,\zeta) + 
					\bar{\mu}_j'(\chi)\bar{N}_{i,j}(\chi,\zeta)$} & = 
				\nonumber \\
				-\int\limits_0^1 \bar{\Theta}_i(\chi,y)\bar{M}_j(\chi,\zeta,y)dy - 
				\sum_{\ell=1}^m \bar{\Psi}_{i,\ell}(\chi) \bar{N}_{\ell,j}(\chi,\zeta), 	
				\label{eq:infmobskalt2}
			\end{align}
		\end{subequations}
		with boundary conditions 
		\begin{subequations}
			\label{eq:infmobskbcalt}
			\begin{align}
				\forall j:& &  \bar{M}_j(\chi,\chi,y) & = \frac{\bar{W}_j(\chi,y)}{\bar{\mu}_j(\chi) + 
					\bar{\lambda}(\chi,y)}, \\
				\forall i \neq j: & & \bar{N}_{i,j}(\chi,\chi) & =\frac{\bar{\Psi}_{i,j}(\chi)}{\bar{\mu}_j(\chi) - 
				\bar{\mu}_i(\chi)}, \\
				\forall i \geq j: & & \bar{N}_{i,j}(\chi,0) & =  \int\limits_0^1 R_i(y)\bar{M}_j(\chi,0,y)dy, \\
				\forall	i < j: & &  \quad \bar{N}_{i,j}(1,\zeta) & = \bar{n}_{i,j}(\zeta),
			\end{align}
		\end{subequations}
		where, by \eqref{eq:nij}, $\bar{n}_{i,j}$ satisfies the compatibility condition
		\begin{equation}
			\bar{n}_{i,j}(1) = \frac{\bar{\Psi}_{i,j}(1)}{\bar{\mu}_j(1) - \bar{\mu}_i(1)}.
		\end{equation}
		
		The equations \eqref{eq:infmobskalt}, \eqref{eq:infmobskbcalt} are  of the same form as the 
		control kernel equations for $\mathbf{K}, \mathbf{L}$ recalled in \ref{app:sf}.
		In particular, the characteristic curves of \eqref{eq:infmobskalt2} are 
		given~by
		\begin{equation}
			\bar{\rho}_j^p(\chi) = \bar{\phi}_p^{-1}(\bar{\phi}_j(\chi)),
		\end{equation}
		for $1 \leq j \leq p \leq m$, where
		\begin{equation}
			\bar{\phi}_j(\chi) = \int\limits_0^\chi \frac{ds}{\bar{\mu}_j(s)}, \qquad j = 1,\ldots,m.
		\end{equation}
		The characteristic curves of \eqref{eq:infmobskalt2} are the potential discontinuity lines of 
		$\bar{N}_{i,j}$ for $i>j$, but we can split the kernels into subdomains of continuity by 
		segmenting 
		the domain $\mathcal{T}$ of the kernel equations as
		\begin{equation}
			\label{eq:Tseg}
			\bar{\mathcal{T}}_j^p = \left\{(\chi,\zeta) \in [0,1]^2: \bar{\rho}_j^{p+1}(\chi) \leq \zeta \leq 
			\bar{\rho}_j^p(\chi) 
			\right\},
		\end{equation}
		for $1 \leq j \leq p \leq m$, where we denote $\bar{\rho}_j^{m+1} = 0$ for all $j=1,\ldots,m$. 
		Now, 
		it follows by \cite[Thm 2]{HumBek24arxivc} that \eqref{eq:infmobskalt}, 
		\eqref{eq:infmobskbcalt} 
		has a well-posed solution  $\bar{\mathbf{M}} \in L^\infty(\mathcal{T}; L^2([0,1]; 
		\mathbb{R}^{1\times m})), 
		\bar{\mathbf{N}} \in L^\infty(\mathcal{T}; \mathbb{R}^{m\times m})$, where 
		$\bar{\mathbf{M}}, 
		\bar{\mathbf{N}}$ 
		are continuous in $(\chi,\zeta) \in \bar{\mathcal{T}}_j^p$ for all $1 \leq j \leq p \leq m$, and 
		hence, 
		piecewise continuous in $(\chi,\zeta) \in \mathcal{T}$. Consequently, $\mathbf{M}, 
		\mathbf{N}$ 
		given by \eqref{eq:MNalt} is the piecewise continuous solution to \eqref{eq:infmobsk}, 
		\eqref{eq:infmobskbc}.
	\end{proof}
\end{theorem}

\subsection{Stability of the Closed-Loop System}

\begin{theorem}
	\label{thm:infmclstab}
	Under Assumption~\ref{ass:infm}, the closed-loop system \eqref{eq:infm}, \eqref{eq:infmbc} 
	under the observer-based output-feedback control law 
	\eqref{eq:infmUobs}--\eqref{eq:infmobsbc} is exponentially stable on $E_c 
	\times E_c$.
\end{theorem}

The proof is presented at the end of this subsection by utilizing the well-posedness of the 
closed-loop system given by Proposition~\ref{prop:infmclwp} in \ref{app:infmclwp}, 
which allows utilization of the separation principle. We first prove the exponential 
stability of the estimation error dynamics.

\begin{lemma}
	\label{lem:infmobstsstab}
	Under Assumption~\ref{ass:infm}, the estimation error dynamics for $\tilde{u} = \hat{u}-u$ and 
	$\tilde{\mathbf{v}} = \hat{\mathbf{v}} - \mathbf{v}$ given by
	\begin{subequations}
		\label{eq:infmerr}%
		\begin{align}
			\tilde{u}_t(t,x,y) + \lambda(x,y)\tilde{u}_x(t,x,y) + \mathbf{P}_+(x,y)\tilde{\mathbf{v}}(t,0) & 
			= 
			\nonumber \\
			\int\limits_0^1\sigma(x,y,\eta)\tilde{u}(t,x,\eta)d\eta +
			\mathbf{W}(x,y)\tilde{\mathbf{v}}(t,x), & \\
			\tilde{\mathbf{v}}_t(t,x) - \pmb{\Lambda}_-(x)\tilde{\mathbf{v}}_x(t,x) + 
			\mathbf{P}_-(x)\tilde{\mathbf{v}}(t,0)
			& =  \nonumber \\
			\int\limits_0^1\pmb{\Theta}(x,y)\tilde{u}(t,x,y)dy + \pmb{\Psi}(x)\tilde{\mathbf{v}}(t,x),
			&
		\end{align}
	\end{subequations}
	with boundary conditions 
	\begin{subequations}
		\label{eq:infmerrbc}%
		\begin{align}
			\tilde{u}(t,0,y) & = 0, \\
			\tilde{\mathbf{v}}(t,1) & = \int\limits_0^1 \mathbf{R}(y)\tilde{u}(t,1,y)dy.
		\end{align}
	\end{subequations}
	are exponentially stable on $E_c$.
	\begin{proof}
		We transform \eqref{eq:infmerr}, \eqref{eq:infmerrbc} into the following target system
		\begin{subequations}
			\label{eq:infmobsts}%
			\begin{align}
				\tilde{\alpha}_t(t,x,y) + \lambda(x,y)\tilde{\alpha}_x(t,x,y) & = 
				\nonumber \\
				\resizebox{.93\columnwidth}{!}{$ \displaystyle 
					\int\limits_0^1\sigma(x,y,\eta)\tilde{\alpha}(t,x,\eta)d\eta + \int\limits_0^1\int\limits_0^x 
					D_+(x,\xi,y,\eta)\tilde{\alpha}(t,\xi,\eta)d\xi d\eta$}, & \\
				\tilde{\pmb{\beta}}_t(t,x) - \pmb{\Lambda}_-(x)\tilde{\pmb{\beta}}_x(t,x)
				& =  \nonumber \\
				\resizebox{.93\columnwidth}{!}{$ \displaystyle 
					\int\limits_0^1\pmb{\Theta}(x,y)\tilde{\alpha}(t,x,y)dy + \int\limits_0^1\int\limits_0^x 
					\mathbf{D}_-(x,\xi,y)\tilde{\alpha}(t,\xi,y)d\xi dy$}, &
			\end{align}
		\end{subequations}
		with boundary conditions 
		\begin{subequations}
			\label{eq:infmobstsbc}%
			\begin{align}
				\tilde{\alpha}(t,0,y) & = 0, \\
				\tilde{\pmb{\beta}}(t,1) & = \int\limits_0^1 \mathbf{R}(y)\tilde{\alpha}(t,1,y)dy - 
				\int\limits_0^1 
				\mathbf{H}(\xi)\tilde{\pmb{\beta}}(t,\xi)d\xi,
			\end{align}
		\end{subequations}
		where $D_+ \in L^\infty(\mathcal{T}; L^2([0,1]^2;\mathbb{R})), \mathbf{D}_- \in 
		L^\infty(\mathcal{T}; L^2([0,1]; \\ \mathbb{R}^m))$ are given by
		\begin{subequations}
			\label{eq:infmD}
			\begin{align}
				\mathbf{D}_-(x,\xi,y) & = -\mathbf{N}(x,\xi)\pmb{\Theta}(\xi,y) \nonumber \\
				& \qquad - \int\limits_\xi^x
				\mathbf{N}(x,s)\mathbf{D}_-(s,\xi,y)ds, 	\label{eq:appinfmD-} \\
				D_+(x,\xi,y,\eta) & = -\mathbf{M}(x,\xi,y)\pmb{\Theta}(\xi,\eta) \nonumber \\
				& \qquad - \int\limits_\xi^x 
				\mathbf{M}(x,s,y)\mathbf{D}_-(s,\xi,\eta)ds,	\label{eq:appinfmD+}
			\end{align}
		\end{subequations}
		and $\mathbf{H} \in L^\infty([0,1]; \mathbb{R}^{m\times 
			m})$ is strictly upper triangular, i.e., $H_{i,j} = 0$ for all $i \geq j$, with its values for $i<j$ 
			to be 
		determined. The transformation is given by
		\begin{subequations}
			\label{eq:infmobsV}
			\begin{align}
				\tilde{u}(t,x,y) & = \tilde{\alpha}(t,x,y) + \int\limits_0^x 
				\mathbf{M}(x,\xi,y)\tilde{\pmb{\beta}}(t,\xi)d\xi, \label{eq:infmobsV1} \\
				\tilde{\mathbf{v}}(t,x) & = \tilde{\pmb{\beta}}(t,x) + \int\limits_0^x 
				\mathbf{N}(x,\xi)\tilde{\pmb{\beta}}(t,\xi)d\xi, \label{eq:infmobsV2}
			\end{align}
		\end{subequations}
		where $\mathbf{M}, \mathbf{N}$ are the observer kernels. In order for \eqref{eq:infmobsV} to 
		transform \eqref{eq:infmerr}, \eqref{eq:infmerrbc} into \eqref{eq:infmobsts}, 
		\eqref{eq:infmobstsbc}, the observer kernels need to satisfy  
		\eqref{eq:infmobsk}, \eqref{eq:infmobskbc} as shown in \ref{app:obsk}, 
		where we also obtain \eqref{eq:infmP} and \eqref{eq:infmD} for $D_+, \mathbf{D}_-$. Note 
		that as $\mathbf{M} \in L^\infty(\mathcal{T}; L^2([0,1];\mathbb{R}^{1\times m})), \mathbf{N} 
		\in L^\infty(\mathcal{T}; \mathbb{R}^{m\times m})$ are well-posed by 
		Theorem~\ref{thm:obskwp}, \eqref{eq:appinfmD-} is a 
		Volterra equation of second kind for $\mathbf{D}_-$ and it has a well-posed solution 
		$\mathbf{D}_- \in L^\infty(\mathcal{T}; L^2([0,1];\mathbb{R}^m))$ by \cite[Lem. 
		7]{HumBek24arxivc}. Thereafter $D_+ \in L^\infty(\mathcal{T}; L^2([0,1]^2; \mathbb{R}))$ is 
		uniquely determined by \eqref{eq:appinfmD+} based on $\mathbf{C}_-, \mathbf{M}$, and
		$\pmb{\Theta}$. Finally, evaluating \eqref{eq:infmobsV} along $x = 1$ and using the 
		boundary conditions \eqref{eq:infmerrbc}, \eqref{eq:infmobstsbc} gives 
		\begin{equation}
			\mathbf{H}(\xi) = \mathbf{N}(1,\xi) - \int\limits_0^1 \mathbf{R}(y)\mathbf{M}(1,\xi,y)dy,
		\end{equation}
		which splits into \eqref{eq:infmobskbc3} and
		\begin{equation}
			\label{eq:Hij}
			H_{i,j}(\xi) = N_{i,j}(1,\xi) - \int\limits_0^1 R_i(y)M_j(1,\xi,y)dy, \quad \forall i < j,
		\end{equation}
		determining the nonzero entries of $\mathbf{H}$.

		We then establish that the transform \eqref{eq:infmobsV} is boundedly invertible, which 
		follows as 
		\eqref{eq:infmobsV2} is a Volterra equation of second kind for $\tilde{\pmb{\beta}}(\cdot,x)$ 
		in 
		terms of $\tilde{\mathbf{v}}(\cdot,x)$ and $\mathbf{N}$. Thus, by \cite[Thm 
		2.3.6]{HocBook}, 
		\eqref{eq:infmobsV2} has a unique solution 
		$\tilde{\pmb{\beta}}(t,\cdot) \in L^2([0,1];\mathbb{R}^m)$ for all $t \geq 0$. Thereafter, 
		$\tilde{\alpha}(\cdot,t) \in L^2([0,1]; L^2([0,1];\mathbb{R}))$ is uniquely determined by 
		\eqref{eq:infmobsV1} for all $t \geq 0$. Thus, the transform \eqref{eq:infmobsV} is boundedly 
		invertible, and hence, the exponential stability of \eqref{eq:infmerr}, \eqref{eq:infmerrbc} is 
		equivalent to the exponential stability of \eqref{eq:infmobsts}, \eqref{eq:infmobstsbc}.
		
		In order to show that \eqref{eq:infmobsts}, \eqref{eq:infmobstsbc} is exponentially stable, 
		consider a Lyapunov functional with parameters $\delta, \mathbf{B} = 
		\operatorname{diag}(B_1, 
		\ldots, B_m) > 0$ of the form 
		\begin{align}
			\label{eq:lyapobs}
			V(t) & = \int\limits_0^1 \int\limits_0^1 e^{-\delta
				x}\frac{\tilde{\alpha}^2(t,x,y)}{\lambda(x,y)}dydx \nonumber \\
			& \quad  + \int\limits_0^1
			e^{\delta 
			x}\tilde{\pmb{\beta}}^T(t,x)\mathbf{B}\pmb{\Lambda}_-^{-1}(x)\tilde{\pmb{\beta}}(t,x)dx.
		\end{align}
		Computing $\dot{V}(t)$ and integrating by parts in $x$ gives 
		\begin{align}
			\label{eq:lyapobsd}
			\dot{V}(t)
			& = \left[-e^{-\delta x}\|\tilde{\alpha}(t,x,\cdot)\|^2_{L^2} +
			e^{\delta x}\|\tilde{\pmb{\beta}}(t,x)\|^2_{\mathbf{B}}
			\right]_0^1 \nonumber 
		\end{align}	
		\begin{align}
			& \quad - \delta\int\limits_0^1 \left(e^{-\delta
				x}\|\tilde{\alpha}(t,x,\cdot)\|_{L^2}^2 + e^{\delta
				x}\|\tilde{\pmb{\beta}}(t,x)\|_{\mathbf{B}}^2 \right)dx \nonumber \\
			& \quad \resizebox{.86\columnwidth}{!}{$ \displaystyle
				+ 2\int\limits_0^1 \int\limits_0^1 \int\limits_0^1e^{-\delta
					x}\frac{\tilde{\alpha}(t,x,y)}{\lambda(x,y)}\sigma(x,\eta,y)\tilde{\alpha}(t,x,\eta)
				d\eta dy dx$} \nonumber
		\end{align}	
		\begin{align}		
			& \quad \resizebox{.7\columnwidth}{!}{$ \displaystyle
				+ 2\int\limits_0^1 \int\limits_0^1 \int\limits_0^1
				\int\limits_0^x e^{-\delta x}
				\frac{\tilde{\alpha}(t,x,y)}{\lambda(x,y)}D_+(x,\xi,y,\eta)\tilde{\alpha}(t,\xi,\eta)$}
			d\xi d\eta dy dx \nonumber \\
			& \quad + \int\limits_0^1 \int\limits_0^1 e^{\delta 
				x} 
			\left(\tilde{\pmb{\beta}}^T(t,x)\mathbf{B}\pmb{\Lambda}_-^{-1}(x)\pmb{\Theta}(x,y) \right. 
			\nonumber \\
			& \left. \qquad  + 
			\pmb{\Theta}^T(x,y)\mathbf{B}\pmb{\Lambda}_-^{-1}(x)\tilde{\pmb{\beta}}(t,x)\right)\tilde{\alpha}(t,x,y)dydx
			\nonumber \\
			& \quad + \int\limits_0^1 \int\limits_0^1 \int\limits_0^x
			e^{\delta 
				x}\left(\tilde{\pmb{\beta}}^T(t,x)\mathbf{B}\pmb{\Lambda}_-^{-1}(x)\mathbf{D}_-(x,\xi,y) 
			\right. \nonumber \\
			& \left. \qquad + \mathbf{D}^T_-(x,\xi,y)\mathbf{B}\pmb{\Lambda}_-^{-1}(x)
			\tilde{\pmb{\beta}}(t,x)\right)\tilde{\alpha}(t,\xi,y)d\xi dy dx,
		\end{align}
		where $\|\cdot\|_{\mathbf{B}}^2 = \left\langle \cdot, \mathbf{B}\cdot
		\right\rangle_{\mathbb{R}^m}$ denotes the $\mathbf{B}$-weighted inner product. Using the
		following bounds (that exist by Assumption~\ref{ass:infm} and Theorem~\ref{thm:obskwp})
		\begin{subequations}
			\label{eq:lyapb}%
			\begin{align}
				m_{\lambda} & = \min_{x,y \in [0,1]} \lambda(x,y), \\
				m_{\mu}  &= \min_{j \in \left\{ 1,\ldots,m \right\}}\min_{x\in
					[0,1]} \mu_j(x), \\
				M_{\sigma} & = \max_{x\in [0,1]} \left\|
				\int\limits_0^1\sigma(x,\cdot,\eta)d\eta
				\right\|_{L^2}, \label{eq:Ms} \\
				M_\Theta & = \max_{j=1,\ldots,m} \max_{x\in[0,1]} \|\Theta_j(x,\cdot)\|_{L^2}, \\
				M_{D_+} & = \esssup_{(x,\xi) \in \mathcal{T}} \left\|
				\int\limits_0^1D_+(x,\xi,\cdot,\eta)d\eta \right\|_{L^2},
				\\
				M_{D_-} & = \max_{j\in \left\{ 1,\ldots,m \right\}}\esssup_{(x,\xi)
					\in \mathcal{T}} \|D_j^-(x,\xi,\cdot)\|_{L^2}, \\
				M_H & = \max_{i,j\in \left\{ 1,\ldots,m \right\}}\esssup_{x\in
					[0,1]}\left| H_{ij}(x) \right|, \\
				M_R & = \max_{j=1,\ldots,m}\|R_j\|_{L^2}, \\
				M_B & = \max_{j=1,\ldots,m} B_j,
			\end{align}
		\end{subequations}
		the boundary conditions \eqref{eq:infmobstsbc}, the Cauchy-Schwartz
		inequality, and $2\left\langle f,g \right\rangle_{L^2} \leq
		\|f\|_{L^2}^2 +\|g\|^2_{L^2}$ for any $f,g \in L^2$, we can estimate \eqref{eq:lyapobsd}
		as 
		\begin{align}
			\label{eq:lyapbest}
			\dot{V}(t) & \leq -\left(e^{-\delta}  - 2e^{\delta}mM_R^2M_B 
			\right)\|\tilde{\alpha}(t,1,\cdot)\|^2_{L^2} \nonumber \\
			& \quad +2e^{\delta}\int\limits_0^1 
			\|\mathbf{H}(x)
			\tilde{\pmb{\beta}}(t,x)\|_\mathbf{B}^2dx \nonumber \\
			& \quad - \delta\int\limits_0^1 \left(e^{-\delta
				x}\|\tilde{\alpha}(t,x,\cdot)\|_{L^2}^2 + e^{\delta
				x}\|\tilde{\pmb{\beta}}(t,x)\|_{\mathbf{B}}^2 \right)dx \nonumber
\end{align}	
\begin{align}				
			& \quad + 2\int\limits_0^1 e^{-\delta 
				x}\frac{M_\sigma+M_{D_+}}{m_\lambda}\|\tilde{\alpha}(t,x,\cdot)\|_{L^2}^2dx \nonumber 
				\\
			& \quad +\displaystyle\int\limits_0^1 me^{\delta 
				x}M_B\frac{M_\Theta^2+M_{D_-}^2}{m_\mu^2}\|\tilde{\alpha}(t,x,\cdot)\|_{L^2}^2dx 
			\nonumber \\
			& \quad + 2\int\limits_0^1 e^{\delta x}\|\tilde{\pmb{\beta}}(t,x)\|_\mathbf{B}^2dx.
		\end{align}
		Due to the triangular structure of $\mathbf{H}$ we can estimate
		\begin{equation}
			\label{eq:HB}
			\|\mathbf{H}(x)\tilde{\pmb{\beta}}(t,x)\|_\mathbf{B}^2 \leq 
			M_H^2\sum_{j=2}^{m}\sum_{\ell=1}^{j-1}  (m-\ell)B_\ell \tilde{\beta}_j^2(t,x),
		\end{equation}
		so that we can enforce 
		\begin{equation}
			2e^\delta \|\mathbf{H}(x)\tilde{\pmb{\beta}}(t,x)\|_\mathbf{B}^2 \leq 
			\|\tilde{\pmb{\beta}}(t,x)\|_\mathbf{B}^2,
		\end{equation}
		by assigning
		\begin{equation}
			\label{eq:Bj}
			B_j = 2e^\delta M_H^2\sum_{\ell=1}^{j-1}(m-\ell)B_\ell, 	\qquad j = 2,\ldots,m.
		\end{equation}
		In order to determine $B_1$, the first term of \eqref{eq:lyapbest} needs to be non-positive, 
		which 
		gives 
		\begin{equation}
			\label{eq:MB}
			M_B \leq \frac{e^{-2\delta}}{2mM_R^2},
		\end{equation}
		where $M_B = B_m$  by \eqref{eq:Bj}, and hence, $B_1$ can be assigned~as
		\begin{equation}
			B_1 = \frac{e^{-2\delta}}{2mM_R^2\sum_{\ell = j}^m \left(2e^\delta M_H^2\right)^{\ell 
					-1}\frac{(m-1)!}{(m-\ell)!}}.
		\end{equation}
		Finally, $\dot{V}(t)$ can be made negative definite by choosing $\delta$ such that
		\begin{equation}
			\delta > \max \left\{2\frac{M_\sigma + M_{D_+}}{m_\lambda} + 
			\frac{M_\Theta^2+M_{D_-}^2}{2M_R^2m_\mu^2},3\right\},
		\end{equation}
		and thus, the dynamics \eqref{eq:infmobsts}, \eqref{eq:infmobstsbc} are exponentially stable.
	\end{proof}
\end{lemma}

\paragraph*{Proof of Theorem~\ref{thm:infmclstab}}
Since the closed-loop system \eqref{eq:infm},~\eqref{eq:infmbc}, 
\eqref{eq:infmUobs}--\eqref{eq:infmobsbc},  is 
well-posed by Proposition~\ref{prop:infmclwp}, we can introduce a change of variables $\tilde{z} =
\hat{z} - z$ and write the closed-loop system equivalently, utilizing the notation of 
\ref{app:infmclwp},~as
\begin{equation}
	\label{eq:infmclerr}  
	\begin{bmatrix}
		\dot{z}(t) \\ \dot{\tilde{z}}(t)
	\end{bmatrix}
	= 
	\begin{bmatrix}
		A_{-1} + B_Q\mathcal{C} + B\mathcal{K} & -B\mathcal{C}_R
		+ B\mathcal{K} \\ 0 & A_{-1} + B\mathcal{C}_R + P\mathcal{C}
	\end{bmatrix}
	\begin{bmatrix}
		z(t) \\ \tilde{z}(t)
	\end{bmatrix},
\end{equation}
where $\dot{z}(t) =
\left(A_{-1}z(t)+B_Q\mathcal{C}+B\mathcal{K}\right)z(t)$
corresponds to the closed-loop dynamics of \eqref{eq:infm},
\eqref{eq:infmbc} under the backstepping state feedback law
$\mathbf{U}(t) = \mathcal{K}z(t) - \mathcal{C}_Rz(t)$ and
$\dot{\tilde{z}}(t)=\left(A_{-1}+B\mathcal{C}_R+P\mathcal{C}\right)\tilde{z}(t)$
corresponds to the estimation error dynamics \eqref{eq:infmerr},
\eqref{eq:infmerrbc}. As these dynamics are exponentially stable by \cite[Thm 
1]{HumBek24arxivc} and Lemma~\ref{lem:infmobstsstab}, respectively, the exponential stability 
of the closed-loop system follows, e.g., by the Gearhart---Pr\"uss---Greiner Theorem \cite[Thm 
V.1.11]{EngNagBook}, and the decay rate is determined by the smaller one of the diagonal 
dynamics.

\section{Output-Feedback Stabilization of Large-Scale $n+m$ Systems 
	Based on Continuum Kernels} \label{sec:obskappr}

\subsection{Large-Scale $n+m$ Systems of Hyperbolic PDEs}

Consider a system of $n+m$ hyperbolic PDEs\footnote{We scale the sums involving the  $n$-part 
states $u^i, i=1,\ldots,n$ by $1/n$ in order 
to make the considerations in the limit $n\to\infty$ more natural, as discussed in \cite[Rem. 
2.2]{HumBek25b}.}
\begin{subequations}
\label{eq:nmm}%
\begin{align}
\mathbf{u}_t(t,x) + \pmb{\Lambda}_+(x)\mathbf{u}_x(t,x)   & 
 	= \frac{1}{n}\pmb{\Sigma}(x)\mathbf{u}(t,x) +
	\mathbf{W}(x)\mathbf{v}(t,x), \\
	\mathbf{v}_t(t,x) - \pmb{\Lambda}_-(x)\mathbf{v}_x(t,x)
	& = 
	\frac{1}{n}\pmb{\Theta}(x)\mathbf{u}(t,x) + \pmb{\Psi}(x)\mathbf{v}(t,x), 
	\end{align}
\end{subequations}
with boundary conditions
\begin{align}
\label{eq:nmmbc}%
\mathbf{u}(t,0) & = \mathbf{Q}\mathbf{v}(t,0), 
& & \mathbf{v}(t,1) = \frac{1}{n}\mathbf{R}\mathbf{u}(t,1)  + \mathbf{U}(t),
\end{align}
where we employ the matrix
notation for $\mathbf{u}, \mathbf{v}, \mathbf{U}, \pmb{\Lambda}_+, \pmb{\Lambda}_-, 
\pmb{\Sigma}$, 
$\mathbf{W}, \pmb{\Theta},
\pmb{\Psi}, \mathbf{Q}$, and $\mathbf{R}$ for the sake of 
conciseness, that is, $\mathbf{u} = \left( u^i \right)_{i=1}^n, \mathbf{v} = \left( v^j \right)_{j=1}^m$,
$\mathbf{U} = \left( U^j \right)_{j=1}^m$, and the parameters are as follows.
\begin{assumption}
\label{ass:nm}
The parameters of \eqref{eq:nmm}, \eqref{eq:nmmbc} are such that
\begin{subequations}
\label{eq:nmparam}
\begin{align}
	\pmb{\Lambda}_+ & = \operatorname{diag}(\lambda_i)_{i=1}^n \in C^1([0,1]; 
	\mathbb{R}^{n\times n)}, \\
	\pmb{\Lambda}_- & = \operatorname{diag}(\mu_j)_{j=1}^m
	\in C^1([0,1]; \mathbb{R}^{m\times m}), \\
	\pmb{\Sigma} & = (\sigma_{i,j})_{i,j=1}^n \in C([0,1];
                       \mathbb{R}^{n\times n}), \\
  \mathbf{W} & = (w_{i,j})_{i=1,}^n{}_{j=1}^m
	 \in C([0,1];\mathbb{R}^{n\times m}), \\
	\pmb{\Theta} & = ( \theta_{j,i} )_{j=1,}^m{}_{i=1}^n \in
	C([0,1];\mathbb{R}^{m\times n}), \\
	\pmb{\Psi} & = \left( \psi_{i,j} \right)_{i,j=1}^m \in C([0,1];
	\mathbb{R}^{m\times m}), \\
	\mathbf{Q} & = (q_{i,j})_{i=1,}^n{}_{j=1}^m
	 \in \mathbb{R}^{n\times m}, \\
	\mathbf{R} & = (r_{j,i})_{j=1,}^m{}_{i=1}^n
	\in \mathbb{R}^{m\times n},
\end{align}
\end{subequations}
where $\lambda_i(x), \mu_j(x) > 0$ for all $x \in [0,1]$ and $i = 1,\ldots,n, j = 1,\ldots,m$. 
Moreover, $\mu_j$ satisfy \eqref{eq:muass} and $\psi_{j,j} = 0$, 
for $j=1,\ldots,m$.
\end{assumption}

\begin{remark}
\label{rem:nmwp}
Under Assumption~\ref{ass:nm}, it can be shown by using the same arguments as in \cite[Prop. 
A.1]{HumBek25b} that the system \eqref{eq:nmm}, \eqref{eq:nmmbc} is well-posed on the 
Hilbert space $E$.
\end{remark}

\subsection{Observer-Based Output-Feedback Controller Based on Continuum Kernels}

The observer-based backstepping output-feedback law to stabilize
\eqref{eq:nmm}, \eqref{eq:nmmbc} based on continuum kernels is of the form
\begin{align}
\label{eq:nmUobs}
\mathbf{U}(t) & = \frac{1}{n}\int\limits_0^1
                \widetilde{\mathbf{K}}(1,\xi)\hat{\mathbf{u}}(t,\xi)d\xi +  
                \int\limits_0^1\widetilde{\mathbf{L}}(1,\xi)\hat{\mathbf{v}}(t,\xi)d\xi \nonumber
  \\
  & \quad -
  \frac{1}{n}\mathbf{R}\hat{\mathbf{u}}(t,1),
\end{align}
where the control gains $\widetilde{\mathbf{K}}, \widetilde{\mathbf{L}}$  are given in 
\eqref{eq:KLappr} and the observer dynamics for $\hat{\mathbf{u}}, \hat{\mathbf{v}}$ are
\begin{subequations}
	\label{eq:nmmobs}%
	\begin{align}
		\hat{\mathbf{u}}_t(t,x) + \pmb{\Lambda}_+(x)\hat{\mathbf{u}}_x(t,x) + 
		\widetilde{\mathbf{P}}_+(x)(\hat{\mathbf{v}}(t,0) - \mathbf{v}(t,0))   & = \nonumber \\
		 \frac{1}{n}\pmb{\Sigma}(x)\hat{\mathbf{u}}(t,x) +
		\mathbf{W}(x)\hat{\mathbf{v}}(t,x), \\
		\hat{\mathbf{v}}_t(t,x) - \pmb{\Lambda}_-(x)\hat{\mathbf{v}}_x(t,x) + 
		\widetilde{\mathbf{P}}_-(x)(\hat{\mathbf{v}}(t,0) - \mathbf{v}(t,0)) 
		& = \nonumber \\
		\frac{1}{n}\pmb{\Theta}(x)\hat{\mathbf{u}} + \pmb{\Psi}(x)\hat{\mathbf{v}}(t,x), 
	\end{align}
\end{subequations}
with boundary conditions
\begin{align}
	\label{eq:nmmobsbc}%
	\hat{\mathbf{u}}(t,0) & = \mathbf{Q}\mathbf{v}(t,0), 
	& & \hat{\mathbf{v}}(t,1) = \frac{1}{n}\mathbf{R}\hat{\mathbf{u}}(t,1)  + \mathbf{U}(t),
\end{align}
where the  output injection gains $\widetilde{\mathbf{P}}_+, 
\widetilde{\mathbf{P}}_-$ are taken based on the continuum observer kernels $\mathbf{M}, 
\mathbf{N}$, satisfying \eqref{eq:infmobsk}--\eqref{eq:nij}, as
\begin{subequations}
	 \label{eq:Pappr}
\begin{align}
\widetilde{\mathbf{P}}_+(x) & = \mathcal{F}_n^*\mathbf{M}(x,0,\cdot)\pmb{\Lambda}_-(0), \\
 \widetilde{\mathbf{P}}_-(x) & = \mathbf{N}(x,0)\pmb{\Lambda}_-(0).
\end{align}
\end{subequations}
We have the following result.

\begin{theorem}
	\label{thm:cak}
	Consider an $n+m$ system \eqref{eq:nmm}, \eqref{eq:nmmbc} with parameters satisfying 
	Assumption~\ref{ass:nm}. Construct respective continuum parameters satisfying 
	Assumption~\ref{ass:infm} and \eqref{eq:nmcap}, and solve the continuum control and 
	observer 
	kernel equations \eqref{eq:infmk}--\eqref{eq:infmkbca1} and 
	\eqref{eq:infmobsk}--\eqref{eq:nij} for $\mathbf{K}, \mathbf{L}$ and $\mathbf{M}, \mathbf{N}$, 
	respectively, under these parameters. When $n$ is sufficiently large, the observer-based 
	output-feedback law \eqref{eq:nmUobs}--\eqref{eq:Pappr} exponentially stabilizes the $n+m$ 
	system \eqref{eq:nmm}, \eqref{eq:nmmbc} on $E$.
\end{theorem}

\subsection{Proof of Theorem~\ref{thm:cak}}

Firstly, we state the following auxiliary result for the estimation error dynamics.
\begin{lemma}
	\label{lem:obsak}
	When $n$ is sufficiently large, the estimation error dynamics for the observer 
	\eqref{eq:nmmobs}--\eqref{eq:Pappr} are exponentially stable on $E$.
	\begin{proof}
	Analogously to the notation of  \ref{app:infmclwp}, denote $\mathbf{z(t)} = (\mathbf{u}(t,\cdot), 
	\mathbf{v}(t,\cdot))$ and define
	\begin{equation}
	\label{eq:Andef}
	\mathbf{Az}(t) = \begin{bmatrix}
		-\pmb{\Lambda}_+ \partial_x & 0 \\ 0 & \pmb{\Lambda}_-\partial_x 
	\end{bmatrix}z(t) + Sz(t),
	\end{equation}
	where $Sz(t)$ corresponds to the right-hand side of \eqref{eq:nmm} and the domain of 
	$\mathbf{A}$ is defined as $\mathcal{D}(\mathbf{A}) = \{\mathbf{z} \in 
	H^1([0,1];\mathbb{R}^{n+m}): \mathbf{u}(0) = 0, \mathbf{v}(1) = 0 \}$,
	and denote by $\mathbf{A}_{-1}$ the unique extension of $\mathbf{A}$ to $E$. Moreover,
 	define output operators $\mathbf{C}\mathbf{z}(t) = \mathbf{v}(0,t), \mathbf{C}_R\mathbf{z}(t) 
 	= \frac{1}{n}	\mathbf{R}\mathbf{u}(t,1)$, and additionally 
	$\mathbf{B} = \begin{bmatrix}
		0 \\ \delta_1\pmb{\Lambda}_-
	\end{bmatrix}, \mathbf{B}_Q  = \begin{bmatrix}
	\delta_0\pmb{\Lambda}_+\mathbf{Q} \\ 0
	\end{bmatrix}$, and $\widetilde{\mathbf{P}} = 
	\mathcal{F}^*P$, where $\mathcal{F}^* = \operatorname{diag}(\mathcal{F}_n^*,I_m)$ with 
	$\mathcal{F}_n^*$ defined in \eqref{eq:Fns} and $P$ is defined in \eqref{eq:PK}.
	 The estimation error dynamics for $\tilde{\mathbf{z}} = \hat{\mathbf{z}} - \mathbf{z}$ can then 
	 be written as
	\begin{align}
		\label{eq:nmest}
		\dot{\tilde{\mathbf{z}}}(t) & = (\mathbf{A}_{-1} + \mathbf{BC}_R + 
		\widetilde{\mathbf{P}}\mathbf{C})\tilde{\mathbf{z}}(t) \nonumber \\
		& = (\mathbf{A}_{-1} + \mathbf{BC}_R + 
		\mathbf{P}\mathbf{C})\tilde{\mathbf{z}}(t) + \Delta \mathbf{PC}\tilde{\mathbf{z}}(t),
	\end{align}
	where we denote $\Delta\mathbf{P} = \mathbf{P} - \widetilde{\mathbf{P}}$, where $\mathbf{P}$ 
	denotes the output injection operator corresponding to the exact $n+m$ observer kernels given 
	in \cite[(60)--(70)]{HuLDiM16}.
	The dynamics \eqref{eq:nmest} are well-posed as $\mathbf{P}, \widetilde{\mathbf{P}}$ are 
	bounded linear operators and 
	$\mathbf{C}$ is an admissible output operator for $\mathbf{A}_P := \mathbf{A}_{-1} + 
	\mathbf{BC}_R + \mathbf{P}\mathbf{C}$ by \cite[Thm 5.4.2]{TucWeiBook}. By duality \cite[Thm 
	4.4.3]{TucWeiBook}, we have that $\mathbf{C}^*$ is an admissible control operator for 
	$\mathbf{A}_P^*$\footnote{Note that $\mathbf{A}^*_P$ is the (unbounded) adjoint of 
	$\mathbf{A}_P$, i.e., they satisfy $\langle \mathbf{A}_P\mathbf{z}_1, 
		\mathbf{z}_2 \rangle_E = \langle \mathbf{z}_1, \mathbf{A}_P^*\mathbf{z}_2 \rangle_E$ with 
		$\mathbf{z}_1 \in \mathcal{D}(\mathbf{A}_P), \mathbf{z}_2 \in 
		\mathcal{D}(\mathbf{A}_P^*)$; whereas $\mathbf{C}^*$ is defined through the dual pairing 
		$\langle \mathbf{Cz}_1, \mathbf{U}\rangle_{\mathbb{R}^m} = \langle \mathbf{z}_1, 
		\mathbf{C}^*\mathbf{U} \rangle_{\mathcal{D}(\mathbf{A}_P), \mathcal{D}(\mathbf{A}_P)^d}$ 
		with $\mathbf{z}_1 \in \mathcal{D}(\mathbf{A}_P)$, where $ \mathcal{D}(\mathbf{A}_P)^d$ 
		denotes the dual space of $\mathcal{D}(\mathbf{A}_P)$ with respect to the pivot space $E$ 
		(see, e.g., \cite[Sect. 2.9--10]{TucWeiBook}).}. Hence, when $\|\Delta 
	\mathbf{P}\|_{\mathcal{L}(\mathbb{R}^m; E)}$ is sufficiently small, we obtain the exponential 
	stability of $(\mathbf{A}_P + \Delta\mathbf{P}\mathbf{C})^* = \mathbf{A}^*_P + 
	\mathbf{C}^*\Delta\mathbf{P}^*$ by \cite[Prop. A.2]{HumBek25b}.
	
	To conclude the proof, we show that $\|\Delta \mathbf{P}\|_{\mathcal{L}(\mathrm{R}^m; E)}$ 
	becomes arbitrarily small when $n$ is sufficiently large. By applying \eqref{eq:MNalt}, we can 
	transform the observer kernel equations 
	\eqref{eq:infmobsk}, \eqref{eq:infmobskbc} for $\mathbf{M}, \mathbf{N}$ into alternative 
	kernel equations \eqref{eq:infmobskalt}, \eqref{eq:infmobskbcalt} and compare them with the 
	respectively transformed $n+m$ observer kernel equations given in 
	\cite[(73), (74)]{HuLDiM16}. As the transformed observer kernel equations are of the 
	same form as the respective control kernel equations, we have by \cite[Lem. 
	5]{HumBek24arxivc} that $\mathcal{F}_n^*\mathbf{M}(x,0,\cdot), \mathbf{N}(x,0)$, for 
	almost all $x \in [0,1]$, tend arbitrarily close to the exact $n+m$ observer kernels evaluated 
	along $\xi = 0$, when $n$ is sufficiently large. Consequently, as $\pmb{\Lambda}_-(0)$ is 
	unaffected by the continuum approximation, $\|\Delta\mathbf{P}\|_{\mathcal{L}(\mathbb{R}^m; 
	E)}$ becomes arbitrarily small when $n$ is sufficiently large. 
	\end{proof}
\end{lemma}

\paragraph*{Proof of Theorem~\ref{thm:cak}}
The remaining steps are analogous to the proof of Theorem~\ref{thm:infmclstab}. That is, the  
well-posedness of the closed-loop system follows by the same steps as in \ref{app:infmclwp} 
after replacing $E_c$ with $E$ and the operators with the ones introduced in the proof of 
Lemma~\ref{lem:obsak}. The resulting (open-loop) transfer function is in fact identical to 
\eqref{eq:oltf}, so that the same well-posedness arguments for the closed-loop system apply as 
at the end of \ref{app:infmclwp}. Hence, as the closed-loop system is well-posed, we can 
introduce $\tilde{\mathbf{z}} = \hat{\mathbf{z}} - \mathbf{z}$ and write the closed-loop system 
equivalently as
\begin{equation}
	\resizebox{!}{.067\columnwidth}{$
	\label{eq:nmclerr}  
	\begin{bmatrix}
		\dot{\mathbf{z}}(t) \\ \dot{\tilde{\mathbf{z}}}(t)
	\end{bmatrix}
	= 
	\begin{bmatrix}
		\mathbf{A}_{-1} + \mathbf{B}_Q\mathbf{C} + \mathbf{B}\widetilde{\mathbf{K}} & 
		-\mathbf{BC}_R
		+ \mathbf{B}\widetilde{\mathbf{K}} \\ 0 & \mathbf{A}_{-1} + \mathbf{BC}_R + 
		\widetilde{\mathbf{P}}\mathbf{C}
	\end{bmatrix}
	\begin{bmatrix}
		\mathbf{z}(t) \\ \tilde{\mathbf{z}}(t)
	\end{bmatrix}$},
\end{equation}
where the diagonal entries are exponentially stable by \cite[Thm 3]{HumBek24arxivc} and 
Lemma~\ref{lem:obsak}, respectively, when $n$ is sufficiently large. Consequently, due to the 
triangular structure, the closed-loop system is exponentially stable when $n$ is sufficiently large.

\section{Output-Feedback Stabilization of Large-Scale $n+m$ Systems Using Continuum 
Observer} \label{sec:obsappr}

In this section, we consider continuum observer-based stabilization of large-scale $n+m$ 
systems. That is, the observer is taken as the continuum observer \eqref{eq:infmobs}, 
\eqref{eq:infmobsbc}, except that the measurement $\mathbf{v}(0,t)$ is taken from the $n+m$ 
system \eqref{eq:nmm}, \eqref{eq:nmmbc}. Moreover, the control law is taken as 
\eqref{eq:infmUobs} based on the continuum observer (and kernels). The motivation is that the 
computation/implementation of
continuum observer-based control law is then independent of $n$, and thus, so is the respective 
computational complexity. This implies that we may 
potentially gain in computational complexity when computing the control law, as opposed to using 
an $n+m$ observer, when $n$ is large, as computational complexity of the $n+m$ observer 
grows with $n$. In particular, the continuum observer-based control law 
does not require reconstruction of the $n+m$ system state. However, the continuum observer 
can provide an (approximate) estimate for the $n+m$ system state by sampling the continuum 
observer state appropriately in $y$, e.g., by applying $\mathcal{F}^*$ to the continuum observer 
state.\footnote{This results in an additional (mean value)
	approximation error that also becomes arbitrarily small when $n$ is sufficiently large (see, e.g., 
	\cite[(C.44),(C.45)]{HumBek25b}): 
	\begin{equation}
		\label{eq:cobserr}
		\|\mathbf{z}(t)-\mathcal{F}^*\hat{z}(t)\|_E \leq \|\mathcal{F}\mathbf{z}(t) - 
		\hat{z}(t)\|_{E_c} + \|\hat{z}(t) - \mathcal{F}\mathcal{F}^*\hat{z}(t)\|_{E_c}.
	\end{equation}
The estimation error in \eqref{eq:cobserr} decays exponentially to zero as it is evident within the 
proof of Theorem~\ref{thm:sgpstab} (see \eqref{eq:clstab}).}

We state next the main result of the section, which is exponential stability of the closed-loop 
system under the continuum observer-based control law. The proof of this result essentially relies 
on the fact that the continuum, observer-based output-feedback controller stabilizes a virtual 
continuum system, together with the fact that the solutions of the closed-loop, virtual continuum 
system approximate to arbitrary accuracy the closed-loop solutions of the $n+m$ system (under 
the same control input), for sufficiently large $n$, when the parameters of the $n+m$ and 
continuum systems are connected via \eqref{eq:nmcap} (and thus, for large $n$, the parameters 
of the $n+m$ system can be approximated by the parameters of the continuum).

\begin{theorem}
\label{thm:sgpstab}
Consider the continuum observer-based control law \eqref{eq:infmUobs}--\eqref{eq:infmP}, 
where the measurement $\mathbf{v}(0,t)$ is taken from the $n+m$ system \eqref{eq:nmm}, 
\eqref{eq:nmmbc} and the parameters of \eqref{eq:infmUobs}--\eqref{eq:infmP} are 
connected to the parameters of \eqref{eq:nmm}, \eqref{eq:nmmbc} via \eqref{eq:nmcap}.
When $n$ is sufficiently large, the closed-loop system comprising \eqref{eq:nmm}, 
\eqref{eq:nmmbc} and \eqref{eq:infmUobs}--\eqref{eq:infmP} is exponentially stable on $E\times 
E_c$.
\begin{proof}
\underline{Step 1} (well-posedness and alternative representation of the closed-loop system):
We begin by establishing the well-posedness of the considered closed-loop system. Using the 
notation introduced in the proof of Lemma~\ref{lem:obsak} and \ref{app:infmclwp}, the 
closed-loop system can be written as
\begin{equation}
	\label{eq:cobsnmcl}
	\begin{bmatrix}
	\dot{\mathbf{z}}(t) \\ \dot{\hat{z}}(t)
\end{bmatrix} = \begin{bmatrix}
	\mathbf{A}_{-1}+\mathbf{BC}_R + \mathbf{B}_Q\mathbf{C}  & 
	\mathbf{B}\mathcal{K} -  \mathbf{B}\mathcal{C}_R \\
	P\mathbf{C} + B_Q\mathbf{C} & A_{-1} - P\mathcal{C} + B\mathcal{K} 	\end{bmatrix} 
\begin{bmatrix}
	\mathbf{z}(t) \\ \hat{z}(t)
\end{bmatrix}.
\end{equation}
The well-posedness of the closed-loop system \eqref{eq:cobsnmcl} is established in 
\ref{app:cobsnmwp}. For the convenience of the subsequent analysis, we transform the $n+m$ 
system state to $E_c$ by applying $\mathcal{F}$ to the dynamics from the left and using 
$\mathcal{F}^*\mathcal{F} = I$. Hence, we define the notations\footnote{Note that 
	we have $B^n = B$ and $\mathcal{C}^n = C$, as those operators 
	only act on the $\mathbf{v}$-part of the state, which is unaffected by the transforms 
	$\mathcal{F}$ and $\mathcal{F}^*$.}
\begin{subequations}
	\begin{align}
		z^n & = \mathcal{F}\mathbf{z}, & & A^n = \mathcal{F}\mathbf{A}\mathcal{F}^*, \\
		B^n & = \mathcal{F}\mathbf{B} = B, & & B_Q^n = \mathcal{F}\mathbf{B}_Q, \\
		\mathcal{C}^n & = \mathbf{C}\mathcal{F}^* = \mathcal{C}, & & \mathcal{C}_R^n = 
		\mathbf{C}_R\mathcal{F}^*,
	\end{align}
\end{subequations}
so that the closed-loop system \eqref{eq:cobsnmcl} becomes
\begin{equation}
	\label{eq:cobsnmcl2}
	\resizebox{!}{.055\columnwidth}{$
		\begin{bmatrix}
			\dot{z}^n(t) \\ \dot{\hat{z}}(t)
		\end{bmatrix} = \begin{bmatrix}
			A^n_{-1}+B^n\mathcal{C}^n_R + B^n_Q\mathcal{C}^n  & 
			B^n\mathcal{K} - B^n\mathcal{C}_R \\
			P\mathcal{C}^n + B_Q\mathcal{C}^n & A_{-1} - P\mathcal{C} + B\mathcal{K} 	
		\end{bmatrix} 
		\begin{bmatrix}
			z^n(t) \\ \hat{z}(t)
		\end{bmatrix}$}.
\end{equation}
We note that because $\mathcal{F}$ is an isometry, the magnitude of the state of 
\eqref{eq:cobsnmcl} equals that of \eqref{eq:cobsnmcl2}, and thus, the stability properties follow 
from one another.

For the stability analysis, we introduce a virtual continuum state $z$ with dynamics 
\eqref{eq:infm}, \eqref{eq:infmbc} that serves as a continuum approximation of the $n+m$ state 
$\mathbf{z}$ (see \cite[Thm 4]{HumBek24arxivc}). That is, the parameters of the continuum 
system are connected to the $n+m$ parameters according to \eqref{eq:nmcap}, the input is the 
same as that of the $n+m$ system, and the initial condition is taken as $z_0 = 
\mathcal{F}\mathbf{z}_0$. As the control input \eqref{eq:infmUobs} is a
well-posed output of the well-posed system \eqref{eq:cobsnmcl}, the input is locally $L^2$, and 
hence, the virtual continuum system has a well-posed solution $z \in C([0,+\infty); E_c)$ for all 
$z_0 \in E_c$. However, in order to retain the continuum approximation accuracy, we, in fact, 
reset 
the virtual continuum state at $t = kT$, for $k \in \mathbb{N}$ and some given $T > 0$, 
as $z_{kT} = \mathcal{F}\mathbf{z}(kT) \in E_c$, which, actually, results in $z \in C([(k-1)T, kT); 
E_c)$ for all $k \in \mathbb{N}$.

Let us denote by $\tilde{z} = \hat{z} - z$ the virtual continuum estimation error. Now, writing $z^n 
= z - (z - z^n)$, we get the following virtual error dynamics based on 
\eqref{eq:cobsnmcl2} and using $\mathcal{C}^n = \mathcal{C}$
\begin{equation}
	\label{eq:zcoe}
	\dot{\tilde{z}}(t) = (A_{-1} + B\mathcal{C}_R + P\mathcal{C})\tilde{z}(t) + 
	(B_Q+P)\mathcal{C}(z(t) - z^n(t)),
\end{equation}
which consists of internally exponentially stable dynamics plus a perturbation term depending on 
the continuum approximation error $z - z^n$. Using the virtual estimation error, the virtual 
continuum dynamics under the control law \eqref{eq:infmUobs} can be written as
\begin{equation}
	\label{eq:zc}
	\dot{z}(t) = (A_{-1} + B_Q\mathcal{C} + B\mathcal{K})z(t) + (B\mathcal{K} 
	- B\mathcal{C}_R)\tilde{z}(t),
\end{equation}
which consists of internally exponentially stable dynamics plus a perturbation term depending on 
the virtual estimation error. For the actual $n+m$ system, the continuum observer-based control 
law contains additional error terms due to the kernel and parameter approximation. We write the 
term $\mathcal{C}_R\hat{z}$ as
\begin{align}
\mathcal{C}_R\hat{z}(t) 
& = \mathcal{C}_R^n\hat{z}(t) + (\mathcal{C}_R - \mathcal{C}_R^n)\hat{z}(t) \nonumber \\
& = \mathcal{C}_R^nz^n(t) + \mathcal{C}_R^n(z(t) - z^n(t)) 
+ 
\mathcal{C}_R^n\tilde{z}(t) \nonumber \\
& \quad + (\mathcal{C}_R - \mathcal{C}_R^n)\hat{z}(t),
\end{align}
which consists of the exact state feedback, a continuum approximation error, a virtual estimation 
error, and parameter approximation error (operating on the observer state). Performing the same 
steps with the term $\mathcal{K}\hat{z}(t)$ in the control law, and by using $\hat{z}=\tilde{z}+z$, 
we get dynamics
\begin{align}
	\label{eq:Fz}
	\dot{z}^n(t) & = (A_{-1}^n + B^n_Q\mathcal{C}^n + 
	B^n\mathcal{K}^n)z^n(t) 
	\nonumber \\
	& \quad + B^n(\mathcal{K}^n - \mathcal{C}^n_R)(z(t) - z^n(t)) \nonumber \\
	& \quad + B^n(\mathcal{K} - \mathcal{C}_R)\tilde{z}(t)  + B^n(\Delta \mathcal{K}^n - \Delta 
	\mathcal{C}^n_R)z(t),
\end{align}
where we additionally denote $\Delta \mathcal{K}^n = \mathcal{K} - \mathcal{K}^n$ and $\Delta 
\mathcal{C}_R^n = \mathcal{C}_R - \mathcal{C}_R^n$.\footnote{Here $\mathcal{K}^n$ denotes 
an approximate, continuum-based 
stabilizing control kernel for the $n+m$ system. For example, it can be viewed as an 
approximation of the 
continuum control kernels, in particular,  for step approximation we have $\mathcal{K}^n = 
\mathcal{F}\mathcal{F}^*\mathcal{K}$. Such approximate, continuum-based control kernels 
remain stabilizing for the $n+m$ system as the difference from the exact control kernel $\Delta 
\mathcal{K}^n$ becomes arbitrarily small when $n$ is sufficiently large (see \cite[Lem. 
5]{HumBek24arxivc},  \cite[Thm 3]{HumBek24arxivc}, \cite[Thm 3.4]{HuLDiM16}). 
\label{fn:Kappr}}

\underline{Step 2} (stability properties of the nominal closed-loop system): Combining 
\eqref{eq:zcoe}, \eqref{eq:zc}, \eqref{eq:Fz}, the closed-loop dynamics  for 
$z^e := (\tilde{z}, z,z^n)$ can be written as
\begin{equation}
	\label{eq:cocl}
	\dot{z}^e(t) = A_sz^e(t) + \begin{bmatrix}
		(B_Q+P)\mathcal{C} \\
		0 \\ B^n(\mathcal{K}^n - \mathcal{C}_R^n)
	\end{bmatrix}(z(t) -z^n(t)),
\end{equation}
where we denote
\begin{equation}
	\label{eq:coclA}
	\resizebox{!}{.066\columnwidth}{$
	A_s = \begin{bmatrix}
		A_{-1} + B\mathcal{C}_R + P\mathcal{C} & 0 & 0 \\
	B(\mathcal{K} - \mathcal{C}_R)	& A_{-1} + B_Q\mathcal{C} + B\mathcal{K} & 0 \\
	B^n(\mathcal{K} - \mathcal{C}_R) & B^n(\Delta \mathcal{K}^n - \Delta \mathcal{C}^n_R) 
	& A_{-1}^n + B^n_Q\mathcal{C}^n + B^n\mathcal{K}^n
	\end{bmatrix}$},
\end{equation}
which is exponentially stable due to the exponential stability 
of the diagonal entries by Lemma~\ref{lem:infmobstsstab}, \cite[Thm 1]{HumBek24arxivc}, and 
Footnote~\ref{fn:Kappr}, respectively.
The perturbation term in \eqref{eq:cocl} can be estimated via the continuum approximation error, 
which, however, is also dependent on the closed-loop state. Hence, as there is no guarantee that 
the operator acting on the 
continuum approximation error is small in norm, the state-dependent perturbation 
may negatively affect the exponential closed-loop stability. However, we establish next that, in 
fact, when the continuum approximation error is sufficiently small (for sufficiently large $n$) the 
perturbation term does not destroy exponential stability of the closed-loop system.

Using variation of parameters, the solution to \eqref{eq:cocl} can be written as 
\begin{equation}
	\label{eq:cocli}
	z^e(t)  = \mathbb{T}_tz_0^e + \int\limits_0^t \mathbb{T}_{(t-s)}\begin{bmatrix}
		(B_Q+P)\mathcal{C} \\
		0 \\ B^n(\mathcal{K}^n - \mathcal{C}_R^n)
	\end{bmatrix}(z(s) - z^n(s))ds,
\end{equation}
where by  $\mathbb{T}_t$ we denote the semigroup generated by the 
lower-triangular operator in  \eqref{eq:coclA}. We define a 
perturbation as 
\begin{equation}
\mathbf{d}(t) = \begin{bmatrix}
	\label{eq:wpert}
	d_1(t) \\ d_2(t)
\end{bmatrix} := \begin{bmatrix}
\mathcal{C} \\ \mathcal{K}^n - \mathcal{C}_R^n
\end{bmatrix}(z(t) - z^n(t)),
\end{equation}
and write $z^e(t)  = \mathbb{T}_tz_0^e + \Phi_{t}\mathbf{d}$,
where the input-to-state map $\Phi_t$, defined by
\begin{equation}
	\Phi_t\mathbf{d} :=  \int\limits_0^t \mathbb{T}_{(t-s)}\begin{bmatrix}
		B_Q+P & 0 \\ 0 & 0 \\
		0 & B^n
	\end{bmatrix}\mathbf{d}(s)ds,
\end{equation}
can be bounded independently of $t$ due to the exponential stability of the
semigroup $\mathbb{T}_t$ \cite[Porp. 4.4.5]{TucWeiBook}. Moreover, as $\mathbb{T}_t$ is 
exponentially stable, the first term in \eqref{eq:cocli} decays exponentially to zero. Hence, for 
every $t \geq 0$, there exist some $M, \omega, M_\Phi > 0$ such that
\begin{equation}
	\label{eq:zest1}
	\|z^e(t)\|_{E_c^3} \leq Me^{-\omega t}\|z_0^e\|_{E_c^3} + M_\Phi\|\mathbf{d}\|_{L^2([0,t]; 
	\mathbb{R}^{2m})}.
\end{equation}

\underline{Step 3} (estimation of perturbation due to continuum approximation): The remaining 
step is to estimate the perturbation \eqref{eq:wpert} in terms of the 
continuum approximation error, where we employ Proposition~\ref{prop:outappr} from 
\ref{app:caout}. Moreover, as the input $\mathbf{U}(t) = (\mathcal{K} - \mathcal{C}_R)\hat{z}(t)$ 
is a well-posed output of the closed-loop system \eqref{eq:cobsnmcl2}, there exists an operator 
$\Psi_t^{\mathrm{cl}}$ such that $\mathbf{U} = \Psi_t^{\mathrm{cl}}\left( \begin{smallmatrix}
	z_0^n \\ \hat{z}_0
\end{smallmatrix}\right)$. Moreover, by using
\begin{equation}
	\label{eq:zzhest}
	\left\| \begin{bmatrix}
		\hat{z}(t) \\ z^n(t)
	\end{bmatrix} \right\|_{E_c^2} = \left\| \begin{bmatrix}
		1 & 1 & 0  \\ 0 & 0 & 1
	\end{bmatrix}z^e(t)\right\|_{E_c^2} \leq 2\|z^e(t)\|_{E_c^3},
\end{equation}
 for any fixed $T > 0$, we have the estimate by Proposition~\ref{prop:outappr}
\begin{equation}
	\label{eq:dest}
	\|\mathbf{d}\|_{L^2([0,T]; \mathbb{R}^{2m})} \leq \delta_1 \|\mathbf{z}_0\|_E + 2\delta_2 
	M_{\Psi_T^\mathrm{cl}} \|z_0^e\|_{E_c^3},
\end{equation}
where $M_{\Psi_T^\mathrm{cl}} = \|\Psi_T^\mathrm{cl}\|_{\mathcal{L}(E_c^2,L^2([0,T]; 
\mathbb{R}^{2m}))}$. Thus, inserting \eqref{eq:dest} to \eqref{eq:zest1} and using $ 
\|\mathbf{z}_0\|_E = \|z_0^n\|_{E_c} \leq  \|z_0^e\|_{E_c^3}$,
we have, for any fixed $T > 0$, for $t \in [0,T)$ that
\begin{equation}
	\label{eq:zest}
	\|z^e(t)\| \leq Me^{-\omega t}\|z_0^e\|_{E_c^3} + 
	M_\Phi \left(\delta_1 + 2\delta_2 
	M_{\Psi_T^\mathrm{cl}}\right)\|z_0^e\|_{E_c^3}.
\end{equation}
Moreover, in order to retain the accuracy of the continuum approximation, we reset the virtual 
continuum state at $t = T$ to $z(T) = z_T = z^n(T^-)$\footnote{Notation $T^- = T - 
\bar{\varepsilon}$, where $\bar{\varepsilon} > 0$ is arbitrarily small, indicates the time instance 
before the 
reset. Note that $z^n$ and $\hat{z}$ are, in fact, continuous in time, so that $z^n(T^-) = z^n(T)$ 
and $\hat{z}(T^-) = \hat{z}(T)$, as the resets concern only the virtual continuum state $z$.}. By 
\cite[Thm 4]{HumBek24arxivc} (analogously to 
Proposition~\ref{prop:outappr}) and using 
$2M_{\Psi_T^\mathrm{cl}}\|z_0^e\|_{E_c^3} $ as an upper bound for the input on $t \in [0,T)$, 
there exist some $\delta_3, \delta_4 > 0$ such that 
\begin{equation}
	\label{eq:caerr}
	\|z(T^-) - z^n(T^-)\|_{E_c} \leq \delta_3\|\mathbf{z}_0\|_E + 2\delta_4
	M_{\Psi_T^\mathrm{cl}}\|z_0^e\|_{E_c^3}.
\end{equation}
Thus, after the reset we have
\begin{align}
	\label{eq:reserr}
	\|z_T^e\|_{E_c^3} & \leq \|z^e(T^-)\|_{E_c^3} 
	+ \left\| 
	\left(\begin{smallmatrix}
		z(T^-) - z^n(T^-) \\ z(T^-) - z^n(T^-) \\ 0
	\end{smallmatrix}\right)\right\|_{E_c^3} \nonumber \\	
	& \leq Me^{-\omega T}\|z^e_0\|_{E_c^3} + M_\Phi \left(\delta_1 + 2\delta_2 
	M_{\Psi_T^\mathrm{cl}}\right)\|z_0^e\|_{E_c^3} \nonumber \\
	& \quad + 2(\delta_3 + 2\delta_4
	M_{\Psi_T^\mathrm{cl}})\|z_0^e\|_{E_c^3} \nonumber \\
	&\leq \left(Me^{-\omega T} + 
	(M_\Phi + 2)(1+2M_{\Psi_T^\mathrm{cl}})\delta\right)\|z_0^e\|_{E_c^3},
\end{align}
where we denote $\displaystyle \delta = 
\max_{i=1,\ldots,4} \{\delta_i\}$\footnote{We note that, for any fixed $T > 0$, $\delta$ depends 
only on the continuum 
approximation error of the parameters of the $n+m$ system.} and use $\|\mathbf{z}_0\|_E \leq 
\|z_0^e\|_{E_c^3}$.

\underline{Step 4} (derivation of stability estimate): The estimates 
\eqref{eq:zest}--\eqref{eq:reserr} apply analogously on any interval $t \in [kT, (k+1)T)$, for any $k 
\in \mathbb{N}$, provided that $z_0^e$ is replaced with the initial condition $z_{kT}^e$ for that 
particular time interval. Hence, for any $t \in [kT, (k+1)T)$, we have
\begin{align}
	\label{eq:zestkT}
	\|z^e(t)\| & \leq Me^{-\omega (t-kT)}\|z_{kT}^e\|_{E_c^3} \nonumber \\
	& \quad + 
	M_\Phi \left(1 + 2 
	M_{\Psi_T^\mathrm{cl}}\right)\delta\|z_{kT}^e\|_{E_c^3},
\end{align}
where we can estimate 
\begin{align}
	\label{eq:reserrkT}
	\|z_{kT}^e\|_{E_c^3} 
	& \leq \left(Me^{-\omega T} + 
	(M_\Phi + 2)(1+2M_{\Psi_T^\mathrm{cl}})\delta\right)^k\|z_0^e\|_{E_c^3}.
\end{align}
Now, if we fix $T$ and $n$ sufficiently large such that $Me^{-\omega T} + 
(M_\Phi + 2)(1+2M_{\Psi_T^\mathrm{cl}})\delta < 1 =: c$, combining \eqref{eq:zestkT}, 
\eqref{eq:reserrkT} gives, for all $t \in [kT, (k+1)T]$ with $k \in \mathbb{N}$, 
\begin{align}
	\label{eq:zestd}
		\|z^e(t)\|_{E_c^3} & \leq c^k\left(Me^{-\omega(t-kT)} + (M_\Phi + 
		2)(1+2M_{\Psi_T^\mathrm{cl}})\delta\right) \nonumber \\
		& \quad \times\|z_0^e\|_{E_c^3}.
\end{align}
Moreover, using $c^{k} \leq e^{\frac{\log(c)}{T}t}$ for all $t \in [kT, (k+1)T]$, we have, for all $t 
\geq 0$,
\begin{align}
	\label{eq:zeest}
	\|z^e(t)\|_{E_c^3} & \leq \left(M + (M_\Phi +
	2)(1+2M_{\Psi_T^\mathrm{cl}})\delta\right)e^{-\bar ct}\|z_0^e\|_{E_c^3},
\end{align}
where $\bar{c} := -\frac{\log(c)}{T} > 0$ since $c < 1$, and thus, the closed-loop system 
\eqref{eq:cocl} is exponentially stable on $E_c^3$. Consequently, considering the initalization of 
the virtual continuum state to $z(0) = \mathcal{F}\mathbf{z}_0$ and using \eqref{eq:zzhest} with 
$\|z^n\|_{E_c} = \|\mathbf{z}\|_E$ together with $\|z_0^e\|_{E_c}\leq 3\left\|\begin{bmatrix}
	\mathbf{z}_0 \\ \hat{z}_0
\end{bmatrix} \right\|_{E\times E_c}$ (that follows from the definition of $z^e$), we have that
\begin{align}
	\label{eq:clstab}
		\left\| \begin{bmatrix}
		 \mathbf{z}(t) \\ \hat{z}(t)
	\end{bmatrix} \right\|_{E\times E_c} & \leq 6\left(M + (M_\Phi +
	2)(1+2M_{\Psi_T^\mathrm{cl}})\delta\right) \nonumber \\
	& \quad \times e^{-\bar{c}t}\left\| \begin{bmatrix}
		\mathbf{z}_0 \\ \hat{z}_0
	\end{bmatrix} \right\|_{E\times E_c},
\end{align}
 which concludes the proof.
\end{proof}
\end{theorem}

\section{Numerical Example and Simulation Results} \label{sec:ex}

Consider an $n+2$ system \eqref{eq:nmm}, \eqref{eq:nmmbc} with the following parameters for 
$i,\ell = 1,\ldots,n$
\begin{subequations}
	\label{eq:expnm}
	\begin{align}
		\lambda_i(x) & = 1, \quad \mu_1(x) = 2, \quad \mu_2(x) = 1, \\
		\sigma_{i,\ell}(x) & = x^3(x+1)\left(\frac{i}{n} - \frac{1}{2}\right)\frac{\ell}{n}\left(\frac{\ell}{n}- 
		1\right), \\
		W_{i,1}(x) & = 3\left(\frac{i}{n} - \frac{1}{2}\right),  \quad W_{i,2}(x) = 2\left(\frac{i}{n} - 
		\frac{1}{2}\right), \\
		\theta_{1,i}(x) & = -\frac{3i}{n}\left(\frac{i}{n} - 1\right), \quad
		\theta_{2,i}(x) = -\frac{2i}{n}\left(\frac{i}{n} - 1\right), \\
		\psi_{l,j}(x) & = 0, \quad l,j \in \{1,2\}, \\
		Q_{i,1} & = 8\left(\frac{i}{n} - \frac{1}{2}\right), \quad
		Q_{i,2} = -8\left(\frac{i}{n}-2\right), \\
		R_{1,i} & = \cos\left(2\pi\frac{i}{n}\right), \quad R_{2,i} = 2\frac{i}{n}\left(\frac{i}{n}+5\right).
	\end{align}
\end{subequations}
Based on numerical experiments, the $n+2$ system with parameters \eqref{eq:expnm} is 
open-loop unstable.
Due to the particular structure of the parameters, respective continuum parameters satisfying 
\eqref{eq:nmcap} can be constructed as
\begin{subequations}
	\label{eq:exp}
	\begin{align}
		\lambda(x,y) & = 1, \quad \mu_1(x) = 2, \quad \mu_2(x) = 1, \\
		\sigma(x,y,\eta) & = x^3(x+1)\left(y - \frac{1}{2}\right)\eta\left(\eta- 1\right), \\
		W_1(x,y) & = 3\left(y - \frac{1}{2}\right),  \quad W_2(x,y) = 2\left(y - \frac{1}{2}\right), \\
		\theta_1(x,y) & = -3y(y-1), \quad
		\theta_2(x,y) = -2y(y-1), \\				
		\psi_{i,j}(x) & = 0, \quad i,j \in \{1,2\} \\
		Q_1(y) & = 8\left(y - \frac{1}{2}\right), \quad
		Q_2(y) = -8(y-2), \\
		R_1(y) & = \cos(2\pi y), \quad R_2(y) = 2y(y+5).
	\end{align}
\end{subequations}
The observer kernel equations \eqref{eq:infmobsk}--\eqref{eq:nij} have explicit solution
\begin{subequations}
	\begin{align}
		M_1^1(x,\xi,y) & = \left(y - \frac{1}{2}\right), \\
		M_1^2(x,\xi,y) & = e^{x-\frac{1}{2}\xi-1}\left(y - \frac{1}{2}\right), \\
		M_2^2(x,\xi,y) & = e^{x-\xi}\left(y - \frac{1}{2}\right), \\
		N_{1,1}^1(x,\xi) & = 0, \quad  N_{1,1}^2(x,\xi) = 0, \\
		N_{2,1}^1(x,\xi) & = 0, \quad N_{2,1}^2(x,\xi) = e^{x-\frac{1}{2}\xi-1}, \\
		N_{1,2}^2(x,\xi) & = 0, \quad N_{2,2}^2(x,\xi) = e^{x-\xi},
	\end{align}
\end{subequations}
where $M_1^\star(\cdot,y), N_{1,1}^\star$, and $N_{2,1}^\star$ are defined on $\mathcal{T}_1^1
= \{(x,\xi) \in \mathcal{T}: \xi \geq 2x-1 \}$ and $\mathcal{T}_1^2
= \{(x,\xi) \in \mathcal{T}: \xi \leq 2x-1 \}$ for the respective superindex $\star = 1,2$, while 
$M_2^2(\cdot,y), N^2_{2,1}$ and $N^2_{2,2}$ are defined on $\mathcal{T}_2^2 = 
\mathcal{T}$, for each $y\in [0,1]$. Note the discontinuity in $N_{2,1}$ along $\xi = 2x-1$. The 
control kernels $\mathbf{K}, \mathbf{L}$ are the same as in \cite[(75)]{HumBek24arxivc}.

For the simulations, we approximate the $n+2$ system \eqref{eq:nmm}, \eqref{eq:nmmbc} and 
the $n+2$ observer \eqref{eq:nmmobs}--\eqref{eq:Pappr} using 
finite-differences with $128$ grid points in $x$. The continuum observer
\eqref{eq:infmobs}--\eqref{eq:infmP} we implement as an $\hat{n}+2$ system, where $\hat{n}$ is 
a degree of freedom for the numerical implementation, in which we employ finite-difference 
approximations for implementing the observer. For illustrating Theorem~\ref{thm:sgpstab}, we 
emulate the continuum by taking $\hat{n} = 
60$ and consider $n < \hat{n}$. For initial conditions, every state component of the $n+m$ 
system is initialized to $u_0(x) = v_0(x) = \frac{1}{2}\sin(2\pi x)$ and the observer is initialized to 
zero.

\subsection{Illustration of Theorem~\ref{thm:cak}}

\begin{figure}[!ht]
	\includegraphics[width=\columnwidth]{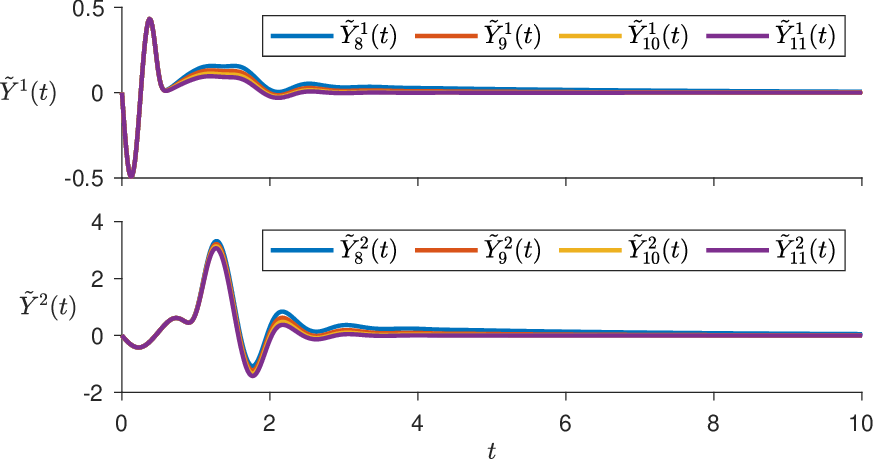}
	\caption{The output estimation errors of the observer \eqref{eq:nmmobs}--\eqref{eq:Pappr}
		for $\mathbf{Y}(t) = \mathbf{v}(0,t)$ from \eqref{eq:nmm}, \eqref{eq:nmmbc} 
		when $n = 8,9,10,11$.}
	\label{fig:Yn}
\end{figure}

\begin{figure}[!ht]
	\includegraphics[width=\columnwidth]{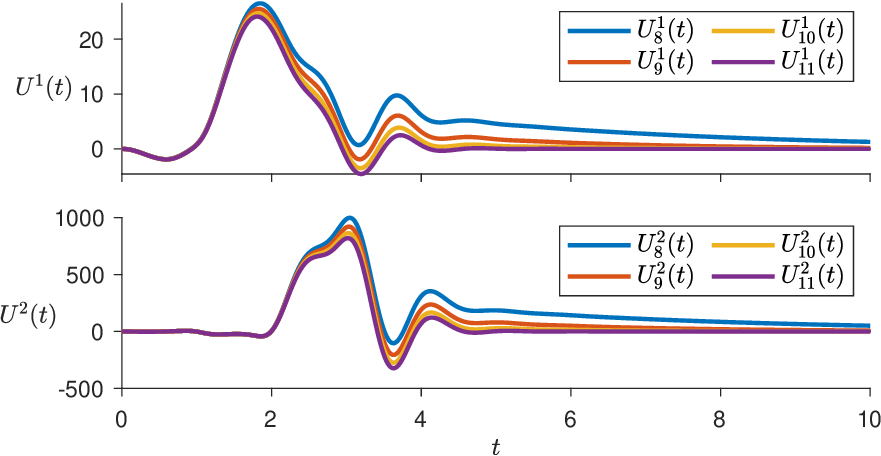}
	\caption{The controls $\mathbf{U}(t)$ based on \eqref{eq:nmUobs}--\eqref{eq:nmmobsbc} in 
	closed-loop with \eqref{eq:nmm}, \eqref{eq:nmmbc} when $n = 8,9,10,11$.}
	\label{fig:Un}
\end{figure}

The simulation results for the output estimation errors and obtained control inputs from 
\eqref{eq:nmUobs}--\eqref{eq:Pappr} are shown in 
Figures~\ref{fig:Yn} and \ref{fig:Un} for $n=8,9,10,11$. Consistently with Theorem~\ref{thm:cak}, 
the estimation errors and the controls tend to zero faster as $n$ increases, as the approximation 
error to the exact kernels becomes smaller as $n$ increases. Moreover, the closed-loop system is 
unstable for $n<7$, which is also consistent with Theorem~\ref{thm:cak}, as large 
approximation errors in the kernels may lead to instability.

\subsection{Illustration of Theorem~\ref{thm:sgpstab}}

\begin{figure}[!ht]
	\includegraphics[width=\columnwidth]{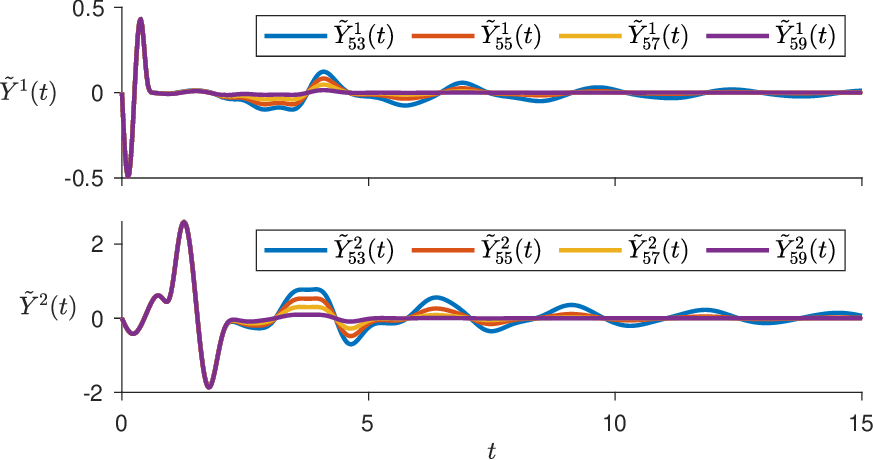}
	\caption{The output estimation errors of the continuum observer 
		\eqref{eq:infmobs}--\eqref{eq:infmP} for $\mathbf{Y}(t) = \mathbf{v}(0,t)$ from 
		\eqref{eq:nmm}, \eqref{eq:nmmbc} 
		when $n = 53,55,57,59$.}
	\label{fig:Y}
\end{figure}

\begin{figure}[!ht]
	\includegraphics[width=\columnwidth]{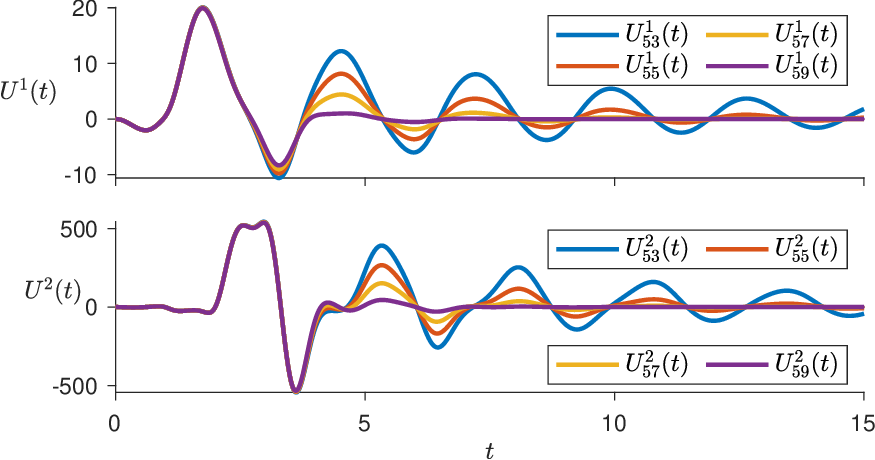}
	\caption{The controls $\mathbf{U}(t)$ based on \eqref{eq:infmUobs}--\eqref{eq:infmobsbc} 
		in closed loop with \eqref{eq:nmm}, \eqref{eq:nmmbc} when $n = 53,55,57,59$.}
	\label{fig:U}
\end{figure}

The simulation results for the output estimation errors and the obtained control inputs from 
\eqref{eq:infmUobs}--\eqref{eq:infmP}, where the continuum is emulated by an $\hat{n}+2$ 
system with $\hat{n}=60$, are shown 
in Figures~\ref{fig:Y} and \ref{fig:U} for $n=53,55,57,59$ in \eqref{eq:nmm}, \eqref{eq:nmmbc}. 
Consistently with Theorem~\ref{thm:sgpstab}, the estimation 
errors and controls tend to zero faster as $n$ tends towards $\hat{n}$. 
Interestingly, the initial transient behavior of the output estimation errors and control 
inputs is virtually the same for all $n$ considered, which may be attributed to the initialization of 
the observer to zero and the transport delays in the dynamics. Finally, we note that, in this case, 
the closed-loop system is unstable for $n < 50$, which is also consistent with 
Theorem~\ref{thm:sgpstab}, as small $n$ may destroy closed-loop stability.

We note that one could also consider the continuum observer implemented as an $\hat{n}+m$ 
system with $\hat{n} < n$, as closed-loop stability is, in fact, retained under any 
approximate observer for which the approximation error related to the solutions  
is sufficiently small, such that $c < 1$ holds in \eqref{eq:zestd}. In particular, implementation of 
the continuum observer as an $\hat{n}+m$ system with $\hat{n} < n$ illustrates a potential gain in 
computational complexity akin to a reduced-order observer. Such a benefit would be more 
evident as $n$ increases, which would allow a continuum observer implementation of order 
$\hat{n}$ (potentially much smaller than $n$) to provide sufficiently accurate estimates of the 
$n+m$ system’s state, without being dominated by errors due to numerical approximation of the 
respective PDEs. Moreover, the accuracy of the continuum approximation depends on the choice 
of the particular numerical scheme employed for solving 
the continuum-based observer. For example, with spectral-based methods (see, e.g., 
\cite{KopBook}) one may obtain results in which the computational complexity benefits are more 
evident, even for smaller $n$. Although numerical computation of the continuum observer 
requires further investigation, the flexibility on the choice of numerical scheme and of the 
observer order $\hat{n}$ may be practically significant, and they are enabled by our 
continuum-based approach for computing control/observer kernels and observer dynamics.

\section{Conclusions and Future Work} \label{sec:conc}

We developed non-collocated, observer-based output-feedback law for a class of 
($\infty+m$) continua of hyperbolic PDE systems. Moreover we employed the developed 
observer (and the respective observer kernels) in designing continuum observer-based 
output-feedback laws for the respective class of large-scale $n+m$ hyperbolic PDEs. The  
motivation of such continuum-based control and observer designs is the potential gain in 
computational 
complexity/flexibility, as the computation of the continuum-based designs is independent of $n$. 
Utilization of the full potential of 
the proposed continuum-based designs calls for development of numerical methods to efficiently 
solve the 
continuum kernel equations \eqref{eq:infmobsk}--\eqref{eq:nij}, 
\eqref{eq:infmk}--\eqref{eq:infmkbca1} and the observer dynamics 
\eqref{eq:infmobs}--\eqref{eq:infmP}. For the former, the power series method 
\cite{VazCheCDC23, HumBek25} could constitute a promising approach, and for the latter, we 
anticipate spectral-based methods (see, e.g., \cite{KopBook}) to be potentially effective. Both of 
these are topics of our ongoing research.

\appendix

\section{Derivation of the Observer Kernels} \label{app:obsk}

Let us first differentiate \eqref{eq:infmobsV} with respect to $x$ and use the Leibniz rule to get 
\begin{align}
	\tilde{u}_x(t,x,y) & = \tilde{\alpha}_x(t,x) + \mathbf{M}(x,x,y)\tilde{\pmb{\beta}}(t,x) \nonumber \\
	& \quad + \int\limits_0^x \mathbf{M}_x(x,\xi,y)\tilde{\pmb{\beta}}(t,\xi)d\xi, \\
	\tilde{\mathbf{v}}_x(t,x) & =  \tilde{\pmb{\beta}}_x(t,x) + 
	\mathbf{N}(x,x)\tilde{\pmb{\beta}}(t,x) \nonumber \\
	& \quad + \int\limits_0^x \mathbf{N}_x(x,\xi)\tilde{\pmb{\beta}}(t,\xi)d\xi.
\end{align}
Moreover, differentiating $\tilde{u}$ in \eqref{eq:infmobsV} with respect to $t$ and using 
\eqref{eq:infmobsts} 
gives
\begin{align}
	\tilde{u}_t(t,x) & = -\lambda(x,y)\tilde{\alpha}_x(t,x,y) + 
	\int\limits_0^1\sigma(x,y,\eta)\tilde{\alpha}(t,x,\eta)d\eta 
	\nonumber \\
	& \qquad + \int\limits_0^x\int\limits_0^1
	D_+(x,\xi,y,\eta)\tilde{\alpha}(t,\xi,\eta)d\eta d\xi \nonumber \\
	& \qquad + \int\limits_0^x 
	\mathbf{M}(x,\xi,y)\pmb{\Lambda}_-(\xi)\tilde{\pmb{\beta}}_\xi(t,\xi)d\xi 
	\nonumber \\
	& \qquad + \int\limits_0^1\int\limits_0^x 
	\mathbf{M}(x,\xi,y)\pmb{\Theta}(\xi,\eta)\tilde{\alpha}(t,\xi,\eta)d\xi d\eta \nonumber \\
	& \qquad  + \resizebox{.73\columnwidth}{!}{$ \displaystyle \int\limits_0^x\int\limits_0^1 
		\int\limits_0^\xi \mathbf{M}(x,\xi,y)
		\mathbf{D}_-(x,s,\eta)\tilde{\alpha}(t,s,\eta)ds d\eta d\xi$},
\end{align}
where integration by parts further gives
\begin{align}
	\int\limits_0^x \mathbf{M}(x,\xi,y)\pmb{\Lambda}_-(\xi)\tilde{\pmb{\beta}}_\xi(t,\xi)d\xi  & = 
	\nonumber \\
	\mathbf{M}(x,x,y)\pmb{\Lambda}_-(x)\tilde{\pmb{\beta}}(t,x) - 	
	\mathbf{M}(x,0,y)\pmb{\Lambda}_-(0)\tilde{\pmb{\beta}}(t,0) & \nonumber \\
	- \int\limits_0^x \left( \mathbf{M}_\xi(x,\xi,y)\pmb{\Lambda}_-(\xi) +  
	\mathbf{M}(x,\xi,y)\pmb{\Lambda}_-'(\xi)\right)\tilde{\pmb{\beta}}(t,\xi) d\xi.
\end{align}
Similarly, differentiating $\tilde{\mathbf{v}}$ in \eqref{eq:infmobsV} with respect to $t$ and using 
\eqref{eq:infmobsts} gives
\begin{align}
	\tilde{\mathbf{v}}_t(t,x) & = \pmb{\Lambda}_-(x)\tilde{\pmb{\beta}}_x(t,x) + 
	\int\limits_0^1\pmb{\Theta}(x,y)\tilde{\alpha}(t,x,y)dy 
	\nonumber \\
	& \qquad + \int\limits_0^x\int\limits_0^1 \mathbf{D}_-(x,\xi,y)\tilde{\alpha}(t,\xi,y)dy d\xi
	\nonumber \\
	& \qquad + \int\limits_0^x \mathbf{N}(x,\xi)\pmb{\Lambda}_-(\xi)\tilde{\pmb{\beta}}_\xi(t,\xi)d\xi 
	\nonumber
\end{align}
\begin{align}		
	& \qquad + \int\limits_0^x\int\limits_0^1 
	\mathbf{N}(x,\xi)\pmb{\Theta}(\xi,y)\tilde{\alpha}(t,\xi,y)dy d\xi \nonumber \\
	& \qquad + \int\limits_0^x\int\limits_0^1 \int\limits_0^\xi
	\mathbf{N}(x,\xi)\mathbf{D}_-(x,s,y)\tilde{a}(t,s,y)ds dy
	d\xi,
\end{align}
where integration by parts further gives
\begin{align}
	\int\limits_0^x \mathbf{N}(x,\xi)\pmb{\Lambda}_-(\xi)\tilde{\pmb{\beta}}_\xi(t,\xi)d\xi & = 
	\nonumber \\
	\mathbf{N}(x,x)\pmb{\Lambda}_-(x)\tilde{\pmb{\beta}}(t,x) - 
	\mathbf{N}(x,0)\pmb{\Lambda}_-(0)\tilde{\pmb{\beta}}(t,0) \nonumber \\
	- \int\limits_0^x  \left(\mathbf{N}_\xi(x,\xi)\pmb{\Lambda}_-(\xi) + 
	\mathbf{N}(x,\xi)\pmb{\Lambda}_-'(\xi)\right)\tilde{\pmb{\beta}}(t,\xi)d\xi.
\end{align}
Thus, in order for \eqref{eq:infmerr} to hold, $\mathbf{M}, \mathbf{N}$ need to satisfy 
\eqref{eq:infmobsk} with boundary conditions \eqref{eq:infmobskbc1}, \eqref{eq:infmobskbc2}, 
along with $\mathbf{P}_+, \mathbf{P}_-$ satisfying \eqref{eq:infmP}, and $D_+, \mathbf{D}_-$ 
satisfying \eqref{eq:infmD}. Moreover, evaluating \eqref{eq:infmobsV} along $x = 1$ and using the 
boundary conditions \eqref{eq:infmerrbc}, \eqref{eq:infmobstsbc} gives  \eqref{eq:infmobskbc3} 
and \eqref{eq:Hij}.

\section{Well-posedness of \eqref{eq:infm}, \eqref{eq:infmbc}, 
	\eqref{eq:infmUobs}--\eqref{eq:infmobsbc}} \label{app:infmclwp}

\begin{proposition}
	\label{prop:infmclwp}
	The closed-loop system \eqref{eq:infm}, \eqref{eq:infmbc}, 
	\eqref{eq:infmUobs}--\eqref{eq:infmobsbc} is well-posed on $E_c \times E_c$.
	\begin{proof}
		In order to write the closed-loop system \eqref{eq:infm}, \eqref{eq:infmbc}, 
		\eqref{eq:infmUobs}--\eqref{eq:infmobsbc} more compactly as an abstract Cauchy problem, 
		define $z(t) = (u(t,\cdot,\cdot), \mathbf{v}(t,\cdot))$, and
		\begin{equation}
			\label{eq:Adef}
			Az(t) = \begin{bmatrix}
				-\lambda \partial_x & 0 \\ 0 & \pmb{\Lambda}_-\partial_x 
			\end{bmatrix}z(t) + Sz(t),
		\end{equation}
		where $Sz(t)$ corresponds to the right-hand side of \eqref{eq:infm} and the domain of $A$ 
		is 
		defined as
		\begin{align}
			\mathcal{D}(A) & = \{z \in H^1([0,1]; L^2([0,1];\mathbb{R}) \times 
			\mathbb{R}^m): \nonumber \\
			& \qquad u(0) = 0, \mathbf{v}(1) = 0 \}.
		\end{align}
		Moreover, denote by $A_{-1}$ the (unique) extension of $A$ from $E_c$ to the dual space of 
		$\mathcal{D}(A)$ with respect to $E_c$, which exists by \cite[Prop. 2.10.2]{TucWeiBook} 
		due to 
		$\mathcal{D}(A)$ being dense in $E_c$. In order to express the boundary couplings, define 
		control
		operators $B, B_Q$ according to the boundary traces in \eqref{eq:infmbc} as \cite[Rem. 
		10.1.6]{TucWeiBook}
		\begin{subequations}
			\begin{equation}
				B  = \begin{bmatrix}
					0 \\ \delta_1\pmb{\Lambda}_-
				\end{bmatrix}, \qquad B_Q  = \begin{bmatrix}
					\delta_0\lambda\mathbf{Q} \\ 0
				\end{bmatrix},
			\end{equation}
		\end{subequations}
		where $\delta_\star$ denotes the Dirac delta function at $x = \star$, and output operators 
		$\mathcal{C}, \mathcal{C}_R$ as
		\begin{equation}
			\mathcal{C}z(t) = \mathbf{v}(0,t), \qquad \mathcal{C}_Rz(t) = \int\limits_0^1 
			\mathbf{R}(y)u(t,1,y)dy.
		\end{equation}
Moreover, define $P, \mathcal{K}$ corresponding to the backstepping observer and controller 
gains as 
\begin{align}
	\label{eq:PK}
	P = \begin{bmatrix}
		\mathbf{P}_+ \\ \mathbf{P}_-
	\end{bmatrix}, & & \mathcal{K}z(t) = \left\langle \begin{bmatrix}
	\mathbf{K}(1,\cdot,\cdot) \\ \mathbf{L}(1,\cdot)
	\end{bmatrix}, \begin{bmatrix}
	u(t) \\ \mathbf{v}(t)
	\end{bmatrix} \right\rangle_{E_c},
\end{align}
so that the closed-loop system  \eqref{eq:infm}, \eqref{eq:infmbc}, 
\eqref{eq:infmUobs}--\eqref{eq:infmobsbc} can be written as an abstract Cauchy problem
\begin{align}
	\label{eq:infmcl}
	\begin{bmatrix}
		\dot{z}(t) \\ \dot{\hat{z}}(t)
	\end{bmatrix} & = \begin{bmatrix}
		A_{-1}+B\mathcal{C}_R + B_Q\mathcal{C}  & -B\mathcal{C}_R + B\mathcal{K} \\
		P\mathcal{C} + B_Q\mathcal{C} & A_{-1} - P\mathcal{C} + B\mathcal{K} 	\end{bmatrix} 
	\begin{bmatrix}
		z(t) \\ \hat{z}(t)
	\end{bmatrix} \nonumber \\
	& = \resizebox{.83\columnwidth}{!}{$\left( \begin{bmatrix}
			A_{-1} & 0 \\ 0 & A_{-1}
		\end{bmatrix} + \begin{bmatrix}
			B & B & B_Q & 0 \\ B & 0 & B_Q & -P
		\end{bmatrix}\begin{bmatrix}
			0 & \mathcal{K} \\
			\mathcal{C}_R & -\mathcal{C}_R \\
			\mathcal{C} & 0 \\ -\mathcal{C} & \mathcal{C}
		\end{bmatrix} \right)\begin{bmatrix}
			z(t) \\ \hat{z}(t)
		\end{bmatrix}$}  \nonumber \\
	& =: \left(A^\mathrm{e}_{-1} + B^\mathrm{e}\mathcal{C}^\mathrm{e} 
	\right)\begin{bmatrix}
		z(t) \\ \hat{z}(t)
	\end{bmatrix}.
\end{align}
Thus, \eqref{eq:infmcl} can be written in an output feedback form, and by \cite[Thm 
13.1.12]{JacZwaBook} the closed-loop system is well-posed if the inverse of $I - 
\mathbf{G}^\mathrm{e}(s)$ exists and is 
bounded for all $\operatorname{Re}(s)$ sufficiently large, where $\mathbf{G}^\mathrm{e}$ is the 
transfer function of the triple $(A^\mathrm{e},B^\mathrm{e},\mathcal{C}^\mathrm{e})$.
By virtue of \cite[Lem. 13.1.14]{JacZwaBook}, the well-posedness of the closed-loop system is 
independent of the (bounded) in-domain coupling terms $Sz(t), S\hat{z}(t)$, meaning that it 
suffices to consider the transfer function $\bar{\mathbf{G}}^\mathrm{e}$, which can be  
computed by \cite[Thm~2.9]{CheMor03} from
\begin{subequations}
	\label{eq:infmcltf}
	\begin{align}
		su(s,x,y) &  = -\lambda(x,y)u_x(s,x,y),  \label{eq:infmcltf1} \\
		s\mathbf{v}(s,x) & = \pmb{\Lambda}_-(x)\mathbf{v}_x(s,x), \label{eq:infmcltf2} \\
		s\hat{u}(s,x,y) &  = -\lambda(x,y)\hat{u}_x(s,x,y) -\mathbf{P}_+(x,y)\mathbf{U}_4(s), \\
		s\hat{\mathbf{v}}(s,x) & = \pmb{\Lambda}_-(x)\hat{\mathbf{v}}_x(s,x) 
		-\mathbf{P}_-(x)\mathbf{U}_4(s), \\
		\mathbf{v}(s,1) & = \mathbf{U}_1(s)+\mathbf{U}_2(s),
		\\
		\hat{\mathbf{v}}(s,1) & = \mathbf{U}_1(s), \\
		u(s,0,y) & = \mathbf{Q}(y)\mathbf{U}_3(s), \\
		\hat{u}(s,0,y) & = \mathbf{Q}(y)\mathbf{U}_3(s), \\
		\mathbf{Y}_1(s) & = \int\limits_0^1\int\limits_0^1
		\mathbf{K}(1,\xi,y)\hat{u}(s,\xi,y)dyd\xi \nonumber \\
		& \qquad + \int\limits_0^1 \mathbf{L}(1,\xi)\hat{\mathbf{v}}(s,\xi) d\xi, \\
		\mathbf{Y}_2(s) & = \int\limits_0^1 \mathbf{R}(y)\left(u(s,1,y)- \hat{u}(s,1,y)\right)dy, \\
		\mathbf{Y}_3(s) & = \mathbf{v}(s,0), \qquad		\mathbf{Y}_4(s) = 		\hat{\mathbf{v}}(s,0) -  
		\mathbf{v}(s,0),
	\end{align}
\end{subequations}
where $\left(\mathbf{Y}_i(s)\right)_{i=1}^4 = 
\bar{\mathbf{G}}^\mathrm{e}(s)\left(\mathbf{U}_i(s)\right)_{i=1}^4$. Since $P$ and $\mathcal{K}$ 
are bounded operators, the components of $\bar{\mathbf{G}}^\mathrm{e}$ involving them 
necessarily tend to zero as $\operatorname{Re}(s)\to\infty$. The remaining components can be 
computed based on the general solution 
\begin{subequations}
	\begin{align}
		u(s,x,y) & = a(y)\exp \left(-s\int\limits_0^x \lambda(\zeta,y)d\zeta\right), \\
		\mathbf{v}(s,x) & = \exp\left(s\int\limits_0^x \pmb{\Lambda}_-(\zeta)d\zeta\right) \mathbf{b},
	\end{align}
\end{subequations}
to \eqref{eq:infmcltf1}, \eqref{eq:infmcltf2} (respectively for $\hat{u}, \hat{\mathbf{v}}$), where 
the coefficients $a, \mathbf{b}$ are solved for from the boundary conditions of 
\eqref{eq:infmcltf}. We obtain (for $\mathbf{U}_4\equiv0$)
\begin{equation}
	\label{eq:oltf}
	\resizebox{.98\columnwidth}{!}{$
		\begin{bmatrix}
			\mathbf{Y}_2(s) \\ \mathbf{Y}_3(s) \\ \mathbf{Y}_4(s) 
		\end{bmatrix} = \begin{bmatrix}
			0 & 0 & 0 \\ 
			\exp \left(-s\int\limits_0^1 
			\pmb{\Lambda}_-(x)dx\right) &  \exp \left(-s\int\limits_0^1 \pmb{\Lambda}_-(x)dx\right) & 
			0 \\
			0 & -\exp \left(-s\int\limits_0^1 
			\pmb{\Lambda}_-(x)dx\right) & 0
		\end{bmatrix}\begin{bmatrix}
			\mathbf{U}_1(s) \\ \mathbf{U}_2(s) \\ \mathbf{U}_3(s)
		\end{bmatrix}$},
\end{equation} 
which concludes that $\operatorname{Re}\bar{\mathbf{G}}_\mathrm{e}(s) \to 0$ as 
$\operatorname{Re}(s) \to \infty$ due to the diagonal matrix $\pmb{\Lambda}_-$ satisfying 
$\pmb{\Lambda}_- > 0$ by Assumption~\ref{ass:infm}. 
Thus, the inverse of $I-\bar{\mathbf{G}}^\mathrm{e}(s)$ exists 
and is bounded when $\operatorname{Re}(s)$ is sufficiently large, and hence, the closed-loop 
system is well-posed by \cite[Thm 13.1.12, Lem. 13.1.14]{JacZwaBook}. That is, for all $z_0, 
\hat{z_0} \in E_c$, there exists a unique solution $z, \hat{z} \in C([0,\infty); E_c)$ to 
\eqref{eq:infmcl}.
\end{proof}	
\end{proposition}

\section{Well-posedness of \eqref{eq:cobsnmcl}} \label{app:cobsnmwp}

The operator in \eqref{eq:cobsnmcl} splits into a feedback form
\begin{align}
\begin{bmatrix}
	\mathbf{A}_{-1}+\mathbf{BC}_R + \mathbf{B}_Q\mathbf{C}  & 
	\mathbf{B}\mathcal{K} -  \mathbf{B}\mathcal{C}_R \\
	P\mathbf{C} + B_Q\mathbf{C} & A_{-1} - P\mathcal{C} + B\mathcal{K} 	
\end{bmatrix} & = \nonumber \\
\begin{bmatrix}
	\mathbf{A}_{-1} & 0	\\
 0 & A_{-1} 
\end{bmatrix} + \begin{bmatrix}
\mathbf{B} & \mathbf{B} & \mathbf{B}_Q & 0 \\ B & 0 & B_Q & -P
\end{bmatrix}\begin{bmatrix}
0 & \mathcal{K} \\
\mathbf{C}_R & -\mathcal{C}_R \\
\mathbf{C} & 0 \\ -\mathbf{C} & \mathcal{C}
\end{bmatrix},
\end{align}
so that the well-posedness of the closed-loop system \eqref{eq:cobsnmcl} can be studied 
similarly to \ref{app:infmclwp}. The changes required to  \ref{app:infmclwp} are, in fact, relatively 
minor and merely involve replacing $u$ with $\mathbf{u}$ and the respective continuum 
operators with the $n+m$ counterparts. For example, \eqref{eq:infmcltf1} becomes
$s\mathbf{u}(s,x) = -\pmb{\Lambda}_+(x)\mathbf{u}_x(s,x)$. By repeating the computations of 
\ref{app:infmclwp}, we get a transfer function that is almost identical to the one shown in 
\eqref{eq:oltf} (ignoring $\mathbf{Y}_1$ and $\mathbf{U}_4$ due to $P$ and $\mathcal{K}$ being 
bounded operators), with the exception of the transfer function from $\mathbf{U}_3$ to 
$\mathbf{Y}_2$ being 
\begin{align}
	\frac{\mathbf{Y}_2(s)}{\mathbf{U}_3(s)} & = 
	\frac{1}{n}\mathbf{R}\exp\left(-s\int\limits_0^1\pmb{\Lambda}_+(\zeta)d\zeta\right)\mathbf{Q} 
	\nonumber \\
	& \quad - \int\limits_0^1 
	\mathbf{R}(y)\exp\left(-s\int\limits_0^1\lambda(\zeta,y)d\zeta\right)\mathbf{Q}(y)dy,
\end{align}
where the real part of both terms tends to zero as $\operatorname{Re}(s) \to \infty$ due to 
$\lambda > 0, \pmb{\Lambda}_+ > 0$ by Assumptions~\ref{ass:infm} and~\ref{ass:nm}. Hence, 
the closed-loop system \eqref{eq:cobsnmcl} is well-posed by the same arguments as those at the 
end of \ref{app:infmclwp}.

\section{Continuum Approximation Result for Well-Posed Outputs} \label{app:caout}

In \cite[Thm 4]{HumBek24arxivc} we have shown that the solution of an $\infty+m$ system 
\eqref{eq:infm}, \eqref{eq:infmbc} can approximate the solution of a respective $n+m$ system 
\eqref{eq:nmm}, \eqref{eq:nmmbc} to arbitrary accuracy on compact time intervals, provided that 
the parameters, initial conditions, and inputs of the two systems are appropriately connected to 
one another. 
In Proposition~\ref{prop:outappr}, we state the respective result for outputs.

\begin{proposition}
	\label{prop:outappr}
	Consider an $n+m$ system \eqref{eq:nmm}, \eqref{eq:nmmbc} with parameters $\mu_j, 
	\psi_{j,\ell}, \theta_{j,i},w_{i,j}, q_{i,j},\lambda_i$, and $\sigma_{i,l}$ for $i,l = 1,\ldots$, $n$ and  
	$j,\ell = 1,\ldots,m$, satisfying Assumption~\ref{ass:nm}, initial conditions $(\mathbf{u}_0, 
	\mathbf{v}_0) \in E$, and a control input $\mathbf{U} \in L_{\rm loc}^2([0,+\infty); 
	\mathbb{R}^m)$. Construct a continuum system \eqref{eq:infm}, \eqref{eq:infmbc} with 
	parameters $\lambda, 
	\mu_j, \sigma, \theta_j, W_j, Q_j, \psi_{j,\ell}$ for $j,\ell=1,\ldots,m$ that satisfy 
	Assumption~\ref{ass:infm} and \eqref{eq:nmcap}, and equip \eqref{eq:infm}, \eqref{eq:infmbc} 
	with initial conditions $(u_0 = \mathcal{F}_n\mathbf{u}_0, \mathbf{v}_0)$ and the same input 
	$\mathbf{U}$. Consider an output space $\mathcal{Y}$ (a Hilbert space) and take an output 
	$\mathbf{Y}^n(t) = \mathcal{C}_\mathrm{a}^n\left(\begin{smallmatrix}
		\mathbf{u}(t) \\ \mathbf{v}(t)
	\end{smallmatrix}\right)$ from the $n+m$ system and the respective output $\mathbf{Y}(t) = 
	\mathcal{C}_\mathrm{a}\left(\begin{smallmatrix}
	u(t)\\ \mathbf{v}(t)
	\end{smallmatrix}\right)$ from the $\infty+m$ system, such that $\mathbf{Y}^n,\mathbf{Y} \in 
	L^2_\mathrm{loc}([0,\infty), \mathcal{Y})$, i.e., such that the outputs are well-posed. On any 
	compact interval $t \in [0,T]$, for any given $T>0$, we have
	\begin{equation}
		\| \mathbf{Y} - \mathbf{Y}^n\|_{L^2([0,T]; \mathcal{Y})} \leq 
		\delta_1\left\|\left(\begin{smallmatrix}
			\mathbf{u}_0\\ \mathbf{v}_0
		\end{smallmatrix}\right)\right\|_E + \delta_2 \|\mathbf{U}\|_{L^2([0,T]; \mathbb{R}^m)},
	\end{equation}
	where $\delta_1, \delta_2 > 0$ become arbitrarily small when $n$ is sufficiently large.
	\begin{proof}
	The proof is analogous to the respective proof for the solutions \cite[Thm 4]{HumBek24arxivc} 
	(see also \cite[Thm 6.1]{HumBek25b}). Namely, considering any well-posed output 
	$\mathbf{Y}$ of the 
	continuum system \eqref{eq:infm}, \eqref{eq:infmbc}, there exist families of bounded linear 
	operators $\Psi_t, \mathbb{F}_t$ for $t \geq 0$, depending continuously on the parameters of 
	\eqref{eq:infm}, \eqref{eq:infmbc}, such that the output is 
	given by
	\begin{equation}
		\label{eq:Yinf}
		\mathbf{Y}  = \Psi_t\left(\begin{smallmatrix}
			u_0\\ \mathbf{v}_0
		\end{smallmatrix}\right)  + \mathbb{F}_t\mathbf{U}.
	\end{equation}
	Now, considering the respective output $\mathbf{Y}^n$ of the $n+m$ system \eqref{eq:nmm}, 
	\eqref{eq:nmmbc} 
	and transforming the system into $E_c$ through $\mathcal{F}$, we get the respective form of 
	the output of the $n+m$ system as
	\begin{equation}
		\label{eq:Yn}
		\mathbf{Y}^n  = \Psi_t^n \left(\begin{smallmatrix}
			\mathcal{F}_n\mathbf{u}_0\\ \mathbf{v}_0
		\end{smallmatrix}\right) + \mathbb{F}_t^n\mathbf{U}.
	\end{equation}
	Subtracting \eqref{eq:Yn} from \eqref{eq:Yinf} and denoting 
$\delta_1 = \|\Psi_T - \Psi_T^n\|_{\mathcal{L}(E_c, L^2([0,T]; \mathcal{Y}))}$, 
		$\delta_2 = \|\mathbb{F}_T - 
		\mathbb{F}_T^n\|_{\mathcal{L}(L^2([0,T]; \mathbb{R}^m), L^2([0,T]; \mathcal{Y}))}$,
	the claim 
	follows after using $u_0 = \mathcal{F}_n\mathbf{u}_0$ and recalling that $\mathcal{F}_n$ is an 
	isometry. In particular, $\delta_1,\delta_2$ become arbitrarily small when $n$ is sufficiently 
	large due to \eqref{eq:nmcap} and the continuous dependence of the output on the parameters.
	\end{proof}
\end{proposition}

\section{Backstepping-Based State-Feedback Stabilzation of \eqref{eq:infm}, 
\eqref{eq:infmbc} and its Application to Stabilization of \eqref{eq:nmm}, \eqref{eq:nmmbc}} 
\label{app:sf}

If the full state information of the system \eqref{eq:infm}, \eqref{eq:infmbc} is available, we have 
shown in \cite[Thm 1]{HumBek24arxivc} that the system can be exponentially  stabilized by the 
state-feedback law
\begin{align}
	\label{eq:infmU}
	\mathbf{U}(t) &  =  \int\limits_0^1\int\limits_0^1
	\mathbf{K}(1,\xi,y)u(t,\xi,y)dyd\xi \nonumber \\
	& \qquad + \int\limits_0^1 \mathbf{L}(1,\xi)\mathbf{v}(t,\xi) d\xi 
	-\int\limits_0^1\mathbf{R}(y)u(t,1,y)dy,
\end{align}
which is of the same form as \eqref{eq:infmUobs} but with the estimated states replaced by the 
actual states. The control kernels $\mathbf{K} \in L^\infty(\mathcal{T}; L^2([0,1]; \mathbb{R}^m)), 
\mathbf{L} \in L^\infty(\mathcal{T}; \mathbb{R}^{m\times m})$ satisfy the control kernel equations
\begin{subequations}
	\label{eq:infmk}%
	\begin{align}
		\pmb{\Lambda}_-(x)\mathbf{K}_x(x,\xi,y) -
		\mathbf{K}_{\xi}(x,\xi,y)\lambda(\xi,y) -
		\mathbf{K}(x,\xi,y)\lambda_{\xi}(\xi,y) & = \nonumber \\
		\mathbf{L}(x,\xi)\pmb{\Theta}(\xi,y) +
		\int\limits_0^1\mathbf{K}(x,\xi,\eta)\sigma(\xi,\eta,y)d\eta,
		& \label{eq:infmkb} \\
		\pmb{\Lambda}_-(x)\mathbf{L}_x(x,\xi) +
		\mathbf{L}_{\xi}(x,\xi)\pmb{\Lambda}_-(\xi) +
		\mathbf{L}(x,\xi)\pmb{\Lambda}_-'(\xi)
		& = \nonumber \\
		\mathbf{L}(x,\xi)\pmb{\Psi}(\xi) + \int\limits_0^1
		\mathbf{K}(x,\xi,y)\mathbf{W}(\xi,y)dy, & \label{eq:infmka}
	\end{align}
\end{subequations}
with boundary conditions
\begin{subequations}
	\label{eq:infmkbc}%
	\begin{align}
	 -\pmb{\Theta}(x,y)  & = \mathbf{K}(x,x,y)\lambda(x,y) + \pmb{\Lambda}_-(x)\mathbf{K}(x,x,y), 
	 \label{eq:kbc1} \\
		\pmb{\Psi}(x) &  =\mathbf{L}(x,x)\pmb{\Lambda}_-(x) - \pmb{\Lambda}_-(x)\mathbf{L}(x,x), 
		\label{eq:kbc2} \\
		L_{i,j}(x,0) & =  \frac{1}{\mu_j(0)} \int\limits_0^1
		K_i(x,0,y)\lambda(0,y)Q_j(y)dy, \quad \forall i \leq j, \label{eq:kbc3} \\
		L_{i,j}(1,\xi) & = l_{i,j}(\xi), \quad \forall j < i, \label{eq:kbca}
	\end{align}
\end{subequations}
where $l_{i,j}$ are arbitrary due to \eqref{eq:kbca} being an artificial boundary condition to 
guarantee well-posedness of the control kernel equations. We 
choose $l_{i,j}$ such that
\begin{equation}
	\label{eq:infmkbca1}%
	l_{i,j}(1) =-\frac{\psi_{i,j}(1)}{\mu_i(1) - \mu_j(1)},
\end{equation}
in order to make the artificial boundary condition compatible with \eqref{eq:kbc2} at $(1,1)$. Note 
that the boundary conditions for $\mathbf{L}(0,0)$ are, in general, overdetermined due to 
\eqref{eq:kbc2} and \eqref{eq:kbc3}, \eqref{eq:kbc1}, which stems potential discontinuities in the 
$\mathbf{L}$ kernels.

Moreover, we have shown in \cite[Sect. 4--5]{HumBek24arxivc} that \eqref{eq:infm}, 
\eqref{eq:infmbc} can be viewed as a continuum
approximation of \eqref{eq:nmm}, \eqref{eq:nmmbc} by connecting the parameters of 
\eqref{eq:infm}, \eqref{eq:infmbc} to the parameters of \eqref{eq:nmm}, 
\eqref{eq:nmmbc} such that $\theta_j, W_j, Q_j,R_j, \lambda$, and $\sigma$ are continuous 
functions satisfying Assumption~\ref{ass:infm} with
\begin{subequations}
	\label{eq:nmcap}%
	\begin{align}
		\theta_j(x,i/n) & = \theta_{j,i}(x), \\
		W_j(x,i/n) & = w_{i,j}(x), \\
		Q_j(i/n) & = q_{i,j}, \\
		R_j(i/n) & = r_{j,i}, \\
		\lambda(x,i/n) & = \lambda_i(x), \\
		\sigma(x,i/n,l/n) & = \sigma_{i,l}(x), \label{eq:nmcape}
	\end{align}
\end{subequations}
for all $x \in [0,1]$, $i,l = 1,\ldots,n$ and $j = 1,\dots,m$.\footnote{As noted in \cite[Sect. 
	IV.A]{HumBek25b}, there are infinitely many ways to construct continuous (or even smooth 
	in $y$) functions satisfying \eqref{eq:nmcap} and Assumption~\ref{ass:infm} based on 
	parameters satisfying Assumption~\ref{ass:nm}.} In particular, by \cite[Thm 
	3]{HumBek24arxivc}, the state-feedback control law
	\begin{align}
		\label{eq:nmU}
		\resizebox{!}{.076\columnwidth}{$\displaystyle 
		\mathbf{U}(t) = \frac{1}{n}\int\limits_0^1
		\widetilde{\mathbf{K}}(1,\xi)\mathbf{u}(t,\xi)d\xi +  
		\int\limits_0^1\widetilde{\mathbf{L}}(1,\xi)\mathbf{v}(t,\xi)d\xi -
		\frac{1}{n}\mathbf{R}\mathbf{u}(t,1)$},
	\end{align}
	exponentially stabilizes the $n+m$ system \eqref{eq:nmm}, \eqref{eq:nmmbc} on $E$,
	where the control gains $\widetilde{\mathbf{K}}, \widetilde{\mathbf{L}}$  are taken based 
	on the continuum control kernels $\mathbf{K}, \mathbf{L}$ in 
	\eqref{eq:infmk}--\eqref{eq:infmkbca1} as
	\begin{align}
		\label{eq:KLappr}
		\widetilde{\mathbf{K}}(x,\xi) & = \mathcal{F}_n^*\mathbf{K}(x,\xi,\cdot), 
		& & \widetilde{\mathbf{L}}(x,\xi) = \mathbf{L}(x,\xi),
	\end{align}
where $\mathcal{F}_n^*$ is given in \eqref{eq:Fns}.


\begin{thebibliography}{10}
	\expandafter\ifx\csname url\endcsname\relax
	\def\url#1{\texttt{#1}}\fi
	\expandafter\ifx\csname urlprefix\endcsname\relax\def\urlprefix{URL }\fi
	\expandafter\ifx\csname href\endcsname\relax
	\def\href#1#2{#2} \def\path#1{#1}\fi
	
	\bibitem{HumBek24arxivc}
	J.-P. {Humaloja}, N.~Bekiaris-Liberis, Backstepping control of continua of
	linear hyperbolic {PDEs} and application to stabilization of large-scale
	$n+m$ coupled hyperbolic {PDE} systems, {arXiv}, 2410.22067 (2024).
	
	\bibitem{HumBek25b}
	J.-P. Humaloja, N.~Bekiaris-Liberis, Stabilization of a class of large-scale
	systems of linear hyperbolic {PDEs} via continuum approximation of exact
	backstepping kernels, IEEE Trans. Automat. Control (2025) 1--16.
	
	\bibitem{BikPhd}
	V.~Bikia, Non-invasive monitoring of key hemodynamical and cardiac parameters
	using physics-based modelling and artificial intelligence, Ph.D. thesis, EPFL
	(2021).
	
	\bibitem{ReyMer09}
	P.~Reymond, F.~Merenda, F.~Perren, D.~Rufenacht, N.~Stergiopulos, Validation of
	a one-dimensional model of the systemic arterial tree, Am. J. Physiol. Heart
	Circ. Physiol. 297 (2009) H208--H222.
	
	\bibitem{ZhaLua22}
	L.~Zhang, H.~Luan, Y.~Lu, C.~Prieur, Boundary feedback stabilization of freeway
	traffic networks: {ISS} control and experiments, IEEE Trans. Control Syst.
	Technol. 30 (2022) 997--1008.
	
	\bibitem{HerKla03}
	M.~Herty, A.~Klar, Modeling, simulation, and optimization of traffic flow
	networks, SIAM J. Sci. Comput. 25 (2003) 1066--1087.
	
	\bibitem{YuHKrs21}
	H.~Yu, M.~Krstic, Output feedback control of two-lane traffic congestion,
	Automatica 125 (2021) Paper No. 109379.
	
	\bibitem{BasCorBook}
	G.~Bastin, J.-M. Coron, Stability and Boundary Sabilization of 1-{D} Hyperbolic
	Systems, Birkh\"{a}user/Springer, [Cham], 2016.
	
	\bibitem{DiaDia17}
	A.~Diagne, M.~Diagne, S.~Tang, M.~Krstic, Backstepping stabilization of the
	linearized {\it {s}aint-{v}enant-{e}xner} model, Automatica 76 (2017)
	345--354.
	
	\bibitem{GuaPri20}
	L.~Guan, C.~Prieur, L.~Zhang, C.~Prieur, D.~Georges, P.~Bellemain, Transport
	effect of {COVID}-19 pandemic in {F}rance, Annu. Rev. Control 50 (2020)
	394--408.
	
	\bibitem{KitBes22}
	C.~Kitsos, G.~Besancon, C.~Prieur, High-gain observer design for a class of
	quasi-linear integro-differential hyperbolic systems-application to an
	epidemic model, IEEE Trans. Automat. Control 67 (2022) 292--303.
	
	\bibitem{HumBek25}
	J.-P. Humaloja, N.~Bekiaris-Liberis, On computation of approximate solutions to
	large-scale backstepping kernel equations via continuum approximation, Syst.
	Control Lett. 196 (2025) 105982.
	
	\bibitem{AllKrs25}
	V.~Alleaume, M.~Krstic, Ensembles of hyperbolic {PDEs}: Stabilization by
	backstepping, IEEE Trans. Automat. Control 70 (2025) 905--920.
	
	\bibitem{DiMVaz13}
	F.~{Di Meglio}, R.~Vazquez, M.~Krstic, Stabilization of a system of {$n+1$}
	coupled first-order hyperbolic linear {PDE}s with a single boundary input,
	IEEE Trans. Automat. Control 58 (2013) 3097--3111.
	
	\bibitem{HuLDiM16}
	L.~Hu, F.~{Di Meglio}, R.~Vazquez, M.~Krstic, Control of homodirectional and
	general heterodirectional linear coupled hyperbolic {PDE}s, IEEE Trans.
	Automat. Control 61 (2016) 3301--3314.
	
	\bibitem{AurBri24}
	J.~Auriol, F.~{Bribiesca Argomedo}, Output-feedback stabilization of $n+m$
	linear hyperbolic {ODE-PDE-ODE} systems, IFAC-PapersOnLine 58 (2024)
	202--207.
	
	\bibitem{Aur24}
	J.~Auriol, Output-feedback stabilization of an underactuated network of {$N$}
	interconnected $n+m$ hyperbolic {PDE} systems, IEEE Trans. Automat. Control
	(2024) 1--12.
	
	\bibitem{EndGab24}
	T.~Enderes, J.~Gabriel, J.~Deutscher, Cooperative output regulation for
	networks of hyperbolic systems using adaptive cooperative observers,
	Automatica 162 (2024) 111506.
	
	\bibitem{AurBre22}
	J.~Auriol, D.~Bresch-Pietri, Robust state-feedback stabilization of an
	underactuated network of interconnected $n+m$ hyperbolic {PDE} systems,
	Automatica 136 (2022) 110040.
	
	\bibitem{AurDiM16}
	J.~Auriol, F.~{Di Meglio}, Minimum time control of heterodirectional linear
	coupled hyperbolic {PDE}s, Automatica 71 (2016) 300--307.
	
	\bibitem{CorHuL17}
	J.-M. Coron, L.~Hu, G.~Olive, Finite-time boundary stabilization of general
	linear hyperbolic balance laws via {F}redholm backstepping transformation,
	Automatica 84 (2017) 95--100.
	
	\bibitem{DiMArg18}
	F.~{Di Meglio}, F.~{Bribiesca Argomedo}, L.~Hu, M.~Krstic, Stabilization of
	coupled linear heterodirectional hyperbolic {PDE}-{ODE} systems, Automatica
	87 (2018) 281--289.
	
	\bibitem{HuLVaz19}
	L.~Hu, R.~Vazquez, F.~{Di Meglio}, M.~Krstic, Boundary exponential
	stabilization of 1-dimensional inhomogeneous quasi-linear hyperbolic systems,
	SIAM J. Control Optim. 57~(2) (2019) 963--998.
	
	\bibitem{BhaShi24}
	L.~Bhan, Y.~Shi, M.~Krstic, Neural operators for bypassing gain and control
	computations in {PDE} backstepping, IEEE Trans. Automat. Control (2024).
	
	\bibitem{WanDia24arxiv}
	S.-S. Wang, M.~Diagne, M.~{Krstic}, Backstepping neural operators for
	$2\times2$ hyperbolic {PDEs}, arXiv, 2312.16762v3 (2024).
	
	\bibitem{QiJZha24}
	J.~Qi, J.~Zhang, M.~Krstic, Neural operators for {PDE} backstepping control of
	first-order hyperbolic {PIDE} with recycle and delay, Syst. Control Lett. 185
	(2024) 105714.
	
	\bibitem{AurMor19}
	J.~Auriol, K.~A. Morris, F.~{Di Meglio}, Late-lumping backstepping control of
	partial differential equations, Automatica 100 (2019) 247--259.
	
	\bibitem{HocBook}
	H.~Hochstadt, Integral {E}quations, {W}iley {C}lassics Edition, John Wiley \&
	Sons, 1989.
	
	\bibitem{EngNagBook}
	K.-J. Engel, R.~Nagel, One-Parameter Semigroups for Linear Evolution Equations,
	Springer, 2000.
	
	\bibitem{TucWeiBook}
	M.~Tucsnak, G.~Weiss, Observation and Control for Operator Semigroups,
	Birkh{\"a}user Verlag AG, 2009.
	
	\bibitem{KopBook}
	D.~Kopriva, Implementing Spectral Methods for Partial Differential Equations,
	Springer Science \& Business Media, 2009.
	
	\bibitem{VazCheCDC23}
	R.~Vazquez, G.~Chen, J.~Qiao, M.~Krstic, The power series method to compute
	backstepping kernel gains: theory and practice, in: IEEE Conf. Decis.
	Control, 2023, pp. 8162--8169.
	
	\bibitem{JacZwaBook}
	B.~Jacob, H.~Zwart, Linear Port-{H}amiltonian Systems on Infinite-dimensional
	Spaces, Birkh{\"a}user, 2012.
	
	\bibitem{CheMor03}
	A.~Cheng, K.~A. Morris, Well-posedness of boundary control systems, SIAM J.
	Control Optim. 42~(4) (2003) 1244--1265.
	
\end{thebibliography}
\end{document}